\documentclass[12pt]{article}
\usepackage{defs1}
\usepackage{slashed}
  
\title{A regulator for smooth manifolds and an index theorem}
\author{Ulrich Bunke\thanks{NWF I - Mathematik,
Universit{\"a}t Regensburg,
93040 Regensburg,
GERMANY, ulrich.bunke@mathematik.uni-regensburg.de}  
}

\theoremstyle{definition}

\newcommand{\incl}{\mathrm{incl}}

\renewcommand{\DD}{\mathbf{DD}}

\newcommand{\bL}{\mathrm{L}}
\newcommand{\bM}{\mathbf{M}}

\newcommand{\bku}{\mathbf{KU}}

\newcommand{\HC}{\mathrm{HC}}

\newcommand{\CC}{\mathbf{CC}}

\renewcommand{\sp}{\mathrm{sp}}
\newcommand{\Iso}{\mathrm{Iso}}

\newcommand{\bs}{\mathbf{s}}

\newcommand{\vol}{{\mathrm{vol}}}

\renewcommand{\Proj}{\mathbf{Proj}}

\newcommand{\cL}{{\mathcal{L}}}

\newcommand{\PSh}{{\mathbf{PSh}}}

\newcommand{\bK}{{\mathbf{K}}}

\newcommand{\Alg}{{\mathbf{Alg}}}

\newcommand{\reg}{{\tt reg}}

\newcommand{\ku}{{\mathbf{ku}}}

\newcommand{\GrCompl}{\mathrm{GrCompl}}
\newcommand{\sto}{\stackrel{\sim}{\to}}

\renewcommand{\Dirac}{\slashed{D}}

\begin{document}
\maketitle
\begin{abstract}
For a  smooth manifold $X$ and  an integer  $d>\dim(X)$ we construct and investigate a natural map
$$\sigma_{d}:K_{d}(C^{\infty}(X))\to \ku\C/\Z^{-d-1}(X)\ .$$ 
Here $K_{d}(C^{\infty}(X))$ is the algebraic $K$-theory group of the algebra of complex valued smooth functions on $X$, and $\ku\C/\Z^{*}$ is the generalized cohomology theory called connective complex $K$-theory with coefficients in $\C/\Z$.

If the manifold $X$ is closed of odd dimension $d-1$ and equipped with a Dirac operator $\Dirac$, then we state and partially prove the conjecture stating that the following two   maps
$K_{d}(C^{\infty}(X))\to \C/\Z$ coincide:
\begin{enumerate}
\item Pair the result of $\sigma_{d}$ with the $K$-homology class of $\Dirac$.
\item Compose the Connes-Karoubi multiplicative character with the classifying map  of the $d$-summable Fredholm module of $\Dirac$.
\end{enumerate}
\end{abstract}
\tableofcontents

\section{Introduction}

The torsion subgroup of the algebraic $K$-theory of the field $\C$ of complex numbers has been calculated
by Suslin \cite{MR772065}. An important tool for this calculation was a 
 collection of homomorphisms 
   \begin{equation}\label{kklajlkjd}r_{2n+1}:K_{2n+1}(\C)\to \C/\Z\end{equation} 
   for $n\in \nat$
which turned out to induce  isomorphisms of torsion subgroups 
$K_{2n+1}(\C)_{tor}\sto \Q/\Z$. 

 \bigskip


%

One may interpret $\C$ as the algebra of complex-valued smooth functions $C^{\infty}(X)$ on the one-point manifold $X=*$.  The first goal of the present paper is to generalize the construction of the homomorphism \eqref{kklajlkjd}  to higher-dimensional smooth manifolds $X$. 
In order to state the result we need the following notation:
\begin{enumerate}
\item
We write $\ku:=\bK^{top}(\C)$ for the connective topological complex $K$-theory spectrum. 
Its homotopy groups are given by
\begin{equation}\label{kckdcjlwdjcljw324234}\pi_{n}(\ku)\cong \left\{\begin{array}{cc} 0 &n<0\\[0.2cm]
\Z&n\ge 0\:\: \mbox{even}\\[0.2cm]
0&n\ge 0\:\: \mbox{odd}\end{array}\right. \ .\end{equation}
\item In general,  we write $\bE^{*}(X)$ for the cohomology groups of the manifold $X$ with coefficients in the spectrum $\bE$.
\item 
If $A$ is an abelian group, then we let $\bM A$ denote the Moore spectrum of $A$ and write
\begin{equation}\label{dewdwedwedewd7879}\bE A:=\bE\wedge \bM A\ .\end{equation}
For example, we can form the spectrum $\ku\C/\Z$. In view of \eqref{kckdcjlwdjcljw324234} and \cite[Eq. (2.1)]{MR551009}
its homotopy groups are given by
\begin{equation}\label{kckdcjlwdjcljw3242341}\pi_{n}(\ku\C/\Z)\cong \left\{\begin{array}{cc} 0 &n<0\\[0.2cm]
\C/\Z&n\ge 0\:\: \mbox{even}\\[0.2cm]
0&n\ge 0\:\: \mbox{odd}\end{array}\right. \ .\end{equation}
\end{enumerate}

\begin{theorem}\label{djekldjlwedwedewdew44dewde}
Let $X$ be a   smooth manifold and $d\in \nat$. If $d>\dim(X)$, then we have a construction of a homomorphism
$$\sigma_{d}:K_{d}(C^{\infty}(X))\to \ku\C/\Z^{-d-1} (X)$$
which is natural in $X$ and induces the map \eqref{kklajlkjd} in homotopy groups for $X=*$
\end{theorem}

\begin{rem}{ \rm
The main point of the theorem is the assertion that there is some interesting generalization of the homomorphism \eqref{kklajlkjd} to higher-dimensional manifolds. In this paper we just give a construction of such a natural homomorphism. We do not address the problem of characterizing it by a collection of natural properties.   Example \ref{ijdoiqwdqwdwqdw} below shows that the map $\sigma_{d}$ contains certain ``higher information''.  As opposed to the case $X=*$, in the higher-dimensional case 
we do not understand its kernel or cokernel.
}
 \end{rem}
 
 \begin{rem}{\rm 
The classical  construction of the homomorphism \eqref{kklajlkjd} relies on the observation by Quillen  that the natural map $K_{*}(\C)\to K_{*}^{top}(\C)$
 from the algebraic to the topological $K$-theory of $\C$
 vanishes rationally in positive degrees. In order to employ this fact for the construction  of $r_{2n+1}$ we work in the stable $\infty$-category $\Sp$ of spectra.  We consider 
the diagram  \begin{equation}\label{djdjk2hdkj2hdkj23hdkj23hdk23jdh23d892736d892}
\xymatrix{\bK(\C)[1..\infty]\ar@{..>}[r]\ar[d]\ar@{-->}[ddr]&\Sigma^{-1}\ku\C/\Z\ar[d]\\\bK(\C) \ar[r]&\ku\ar[d]\\ 
&\ku \C}\ .\end{equation}
Here $\bK(\C)[1..\infty]\to \bK(\C)$ denotes the connective covering of the algebraic $K$-theory spectrum $\bK(\C)$. 
The right column is a segment of a Bockstein fibre sequence for the topological $K$-theory spectrum $ \ku$. Finally,   the middle horizontal map is the canonical map from algebraic to topological $K$-theory.  Now, by Quillen's observation, the  
dashed arrow induces the trivial map in homotopy groups. A choice of a  zero homotopy of this arrow  induces the dotted arrow which in turn  induces the map $r_{2n+1}$ after applying $\pi_{2n+1}(-)$ and identifying $\pi_{2n+1}  (\Sigma^{-1}\ku\C/\Z)$ with $\C/\Z$ using \eqref{kckdcjlwdjcljw3242341}.
Note that the restriction of $r_{2n+1}$ to the torsion subgroup does not depend on the choice of the zero homotopy.
A reference for this construction of \eqref{kklajlkjd}  formulated in a slightly different language  is \cite{MR750694}.
In the present paper we will generalize the alternative  construction \cite[Ex. 6.9]{2013arXiv1311.3188B}.
}\hB \end{rem}

%

 \begin{ex}\label{ijdoiqwdqwdwqdw}{\rm
 For a  unital algebra $A$ let \begin{equation}\label{iotaa}\iota:A^{\times}\to K_{1}(A)\end{equation} be the natural homomorphism
from the  group $A^{\times}$ of units of $A$ to the first algebraic $K$-theory group of $A$.

\bigskip

A complex-valued smooth function $f\in C^{\infty}(S^{1})$  gives rise to an invertible function $\exp(f)\in C^{\infty}(S^{1})$. If 
 $u\in C^{\infty}(S^{1})$ is a second invertible function, then
 we  have two algebraic $K$-theory classes $\iota(\exp(f)),\:\iota(u)\in K_{1}(C^{\infty}(S^{1}))$.
Using the multiplicative structure of the algebraic $K$-theory for commutative algebras we define the class
$$\iota(\exp(f))\cup \iota(u)\in K_{2}(C^{\infty}(S^{1}))\ .$$
 We use the identification  $$\ku\C/\Z^{-3}(S^{1})\cong \ku\C/\Z^{-4}(*)\cong \C/\Z$$ given  by suspension and \eqref{kckdcjlwdjcljw3242341}. With this identification we have  \begin{equation}\label{efewfewfewfewfewfwf234324}\sigma_{2}(\iota(\exp(f))\cup \iota(u))=\left[  \frac{1}{(2\pi i)^{2}}\int_{S^{1}} f\: \frac{du}{u}\right]_{\C/\Z}\ .\end{equation}
 The formula \eqref{efewfewfewfewfewfwf234324}  is a special case of \eqref{hkhkjwdhkwjqdhjkqwdwqdqwd}.
}\hB\end{ex}

\bigskip

In \cite[Example 6.9]{2013arXiv1311.3188B} we explained how one can construct the map \eqref{kklajlkjd}  using techniques of differential cohomology. The 
 differential cohomology approach in particular provides a canonical choice for the dotted arrow in \eqref{djdjk2hdkj2hdkj23hdkj23hdk23jdh23d892736d892}.
The main idea for the construction of $\sigma_{d}$ is to apply the framework of differential cohomology to the algebraic $K$-theory  of the Fr\'echet algebra $C^{\infty}(X)$. The details will be worked out in Section \ref{wdjqdjkwqdnkwqdwqd333}. The final construction of $\sigma_{d}$ will be given in Definition \ref{kdmlkqmdlqwdqwdwqd}.

\bigskip

We now come to the second theme of this paper.
Let us assume that $X$ is a closed  manifold of odd dimension $d$  which carries a generalized Dirac operator $\Dirac$.
This Dirac operator gives rise to a $K$-homology class $[\Dirac]\in \bku_{d}(X)$ and a $d+1$-summable Fredholm module which we will describe in  Subsection \ref{dhwiqdjoqwidqwdqwdqwdwqdqdqwd123} in greater detail. This Fredholm module is classified by a homomorphism $$b_{\Dirac}:C^{\infty}(X)\to \cM_{d}\ ,$$ which is unique up to unitary equivalence, and
where $ \cM_{d}$ denotes the classifying algebra for $d+1$-summable Fredholm modules 
introduced in \cite{MR972606}, see Remark \ref{kljkljddlqwdqwdqwd} for an explicit description. In \cite{MR972606} Connes and Karoubi further introduced 
the multiplicative character $$\delta:K_{d+1}(\cM_{d})\to \C/\Z\ .$$
Since $\ku$ is a $\bku$-module spectrum we can define a map
$$r_{\Dirac}:\ku\C/\Z^{-d-2}(X)\stackrel{\langle-,[\Dirac]\rangle}{\longrightarrow}  \ku\C/\Z^{-2d-2}(*)\stackrel{\eqref{kckdcjlwdjcljw3242341}}{\cong} \C/\Z$$
 given by the pairing with the $K$-homology class $[\Dirac]$.   An explicit construction   of this map using elements of local index theory
  will be given in \eqref{wdqwdwqdwqdwqdq32432432434}.  We now make the following conjecture:
\begin{con}\label{ijdlkqwdwqdqwdqwd}
Assume that $X$ is a closed odd-dimensional manifold of dimension $d$ with a Dirac operator $\Dirac$. Then
the following diagram commutes:
$$\xymatrix{&\ku\C/\Z^{-d-2}(X)\ar[dr]^{ r_{\Dirac}}&\\K_{d+1}(C^{\infty}(X))\ar[ur]^{\sigma_{d+1}}\ar[dr]^{b_{\Dirac}}&&\C/\Z\\
&K_{d+1}(\cM_{d})\ar[ur]^{\delta}& }\ .$$
\end{con}

This conjecture is supported by our  second main result  which asserts that it  holds true if one replaces $K_{d}(C^{\infty}(X))$ by its subgroup of classes which are topologically trivial. 
Note that we do not know any example of a topologically non-trivial class in $K_{n}(C^{\infty}(X))$
for $n> \dim(X)$, see Remark \ref{jwhdqkwdhqwd89qwd79w8q7d}.  In order to explain what topologically trivial means 
 we consider the  homotopification  fibre sequence in spectra (see \eqref{lkdjlkedjwedwed} for details)
$$\bK^{rel}(C^{\infty}(X))\stackrel{\partial}{\to} \bK (C^{\infty}(X))\to \bK^{top}(C^{\infty}(X))\to \Sigma \bK^{rel}(C^{\infty}(X))$$
relating the algebraic $K$-theory spectrum  of $C^{\infty}(X)$ with its topological and relative $K$-theory spectra.
A class in $K_{d+1}(C^{\infty}(X))$ is called topologically trivial if its image in $K_{d+1}^{top}(C^{\infty}(X))$ vanishes, or equivalently, if this class belongs to the image of $\partial:K^{rel}_{d+1}(C^{\infty}(X))\to K_{d+1}(C^{\infty}(X))$. We have the following theorem:
\begin{theorem}\label{dlkdjqwkdjwqlkdqwdqwd}  Assume that $X$ is a closed odd-dimensional manifold of dimension $d$ with a Dirac operator $\Dirac$.  Then the following diagram commutes:
$$\xymatrix{&\ku\C/\Z^{-d-2}(X)\ar[dr]^{ r_{\Dirac}}&\\K^{rel}_{d+1}(C^{\infty}(X))\ar[ur]^{\sigma_{d+1}\circ \partial}\ar[dr]^{b_{\Dirac}\circ \partial}&&\C/\Z\\
&K_{d+1}(\cM_{d})\ar[ur]^{\delta}& }\ .$$
\end{theorem}

\bigskip

In the remainder of this introduction we describe how our constructions are  related with other  results in the literature  relating index and spectral  theory of operators with algebraic $K$-theory of smooth functions.

\bigskip 

Since $C^{\infty}(X)$ is a commutative algebra, the algebraic $K$-theory $K_{*}(C^{\infty}(X))$ is a graded commutative ring. As in Example \ref{ijdoiqwdqwdwqdw}  we can use the map      $$\iota:C^{\infty}(X)^{\times}\stackrel{\eqref{iotaa}}{\to} K_{1}(C^{\infty}(X))$$  
 and the $\cup$-product  in algebraic $K$-theory in order to  construct higher algebraic $K$-theory classes.

 \begin{ex}{\rm  Assume that $X$ is a closed odd-dimensional manifold   with a Dirac operator $\Dirac$.
If $f\in C^{\infty}(X)$, then we can form the unit $e^{f}\in C^{\infty}(X)^{\times}$, and we can consider the class
$\iota(e^{f})\in K_{1}(C^{\infty}(X))$. Note that this element is topologically trivial. Given a collection $f_{1},\dots,f_{d}$ of such smooth functions  we can form the topologically trivial algebraic $K$-theory class
$$\{e^{f_{1}},\dots,e^{f_{d}}\}\in \iota(e^{f_{1}})\cup\dots \cup\iota(e^{f_{d}})\in K_{d}(C^{\infty}(X))\ .$$

The main result of \cite{MR2817643} is an explicit formula \cite[(1.2)]{MR2817643} for the number
$$(\delta\circ  b_{\Dirac})( \{e^{f_{1}},\dots,e^{f_{d}}\}) \in \C/\Z\ .$$
It involves the traces of algebraic expressions build from  the $f_{i}$ and the positive spectral projection of $\Dirac$.
 Kaad's formula can be considered as the analytical side of an index formula. One can interpret our   Conjecture 
   \ref{ijdlkqwdwqdqwdqwd} as providing the topological counterpart. 
   Indeed, since $\{e^{f_{1}},\dots,e^{f_{d}}\}$ is topologically trivial, by Theorem \ref{dlkdjqwkdjwqlkdqwdqwd} we have the equality
      $$(\delta\circ  b_{\Dirac})( \{e^{f_{1}},\dots,e^{f_{d}}\})=(\rho_{\Dirac}\circ \sigma_{d})(\{e^{f_{1}},\dots,e^{f_{d}}\})$$ where  $d-1=\dim(X)$.   
}\end{ex}

 \begin{ex}\label{fijfiowefwefwefewfewfwefewf}{\rm
 We consider the case $X=S^{1}$ with the Dirac operator $\Dirac:=i\partial_{t}$ acting as an unbounded essentially selfadjoint operator with domain $C^{\infty}(S^{1})$ on the Hilbert space $L^{2}(S^{1})$. Let $u_{1},u_{2}\in C^{\infty}(S^{1})^{\times}$ be two  invertible complex-valued functions. Then we form the algebraic $K$-theory class $$\{u_{1},u_{2}\}\in K_{2}(C^{\infty}(S^{1}))\ .$$
Let $P\in B(L^{2}(S^{1}))$ be the projection onto the subspace of positive Fourier modes, i.e. the positive spectral projection of $\Dirac$.
For $f\in C^{\infty}(S^{1})$ we consider the Toeplitz operator $$T_{f}:=PfP\in B(L^{2}(S^{1}))\ ,$$  where $f$ acts as multiplication operator. For two functions $f_{1},f_{2}\in C^{\infty}(S^{1})$ the difference   $T_{f_{1}}T_{f_{2}}-T_{f_{1}f_{2}}$ is a trace class operator.

\bigskip

We let
$\cA\subset B(L^{2}(S^{1}))$ be the algebra generated by all Toeplitz operators $T_{f}$ for $f\in C^{\infty}(S^{1})$ and the algebra of trace class operators $\cL^{1}:=\cL^{1}(L^{2}(S^{1}))$.
We then get the Toeplitz extension
\begin{equation}\label{mdwqkjdqkwdwqdwqdqw}0\to \cL^{1}\to \cA\to C^{\infty}(S^{1})\to 0\ .\end{equation}
Associated to an extension of the trace class operators one has the determinant invariant (see e.g. \cite{MR0390830}) 
$$d:=\det\circ \partial:K_{2}(C^{\infty}(S^{1}))\to \C^{*}\ ,$$
where $\partial:K_{2}(C^{\infty}(S^{1}))\to K_{1}(\cA,\cL^{1})$ is the boundary operator in algebraic $K$-theory associated to the sequence \eqref{mdwqkjdqkwdwqdwqdqw}
and
$\det:K_{1}(\cA,\cL^{1})\to \C^{*}$ is induced by the Fredholm determinant.
The diagram \cite[(3)]{MR2826283} states that
 \begin{equation}\label{wefefwefewfewfewfewfewf}d(\{u_{1},u_{2}\})=\exp(2\pi i \:\delta(b_{\Dirac}(\{u_{1},u_{2}\})) )\ .\end{equation} The determinant invariant was identified by Carey-Pincus with the joint torsion  $\tau(A,B)\in \C^{*}$ which is defined for 
 the pair $A,B$ of Fredolm operators    which commute up to trace class operators.
In the special case of the pair $T_{u_{1}},T_{u_{2}}$ on $\im(P)$ we thus have
   $$d(\{u_{1},u_{2}\})=\tau(T_{u_{1}},T_{u_{2}}) \ .$$  We refer to \cite{2014arXiv1403.4882M} and \cite{MR2854180} for a gentle introduction to joint torsion.

The joint torsion $\tau(T_{u_{1}},T_{u_{2}})$
in turn has been calculated explicitly. 
In the special case
where $u_{1}=e^{f_{1}}$ we have
$$\tau(T_{u_{1}},T_{u_{2}})=\exp(\frac{1}{2\pi i}\int f_{1} \:d\log u_{2})\ .$$
In order to state the result of the calculation in the general case in a comprehensive way we will use the cup product in Deligne cohomology.
Using the isomorphism \eqref{ghjghjdgwqduziqwdwqdwqdqwdwqd} (to be explained below) the invertible functions $u_{i}$ can be interpreted as classes in Deligne cohomology
$H^{1}_{Del}(S^{1},\Z)$. Their cup product is the class
$$u_{1}\cup u_{2}\in H^{2}_{Del}(S^{1},\Z) \ .$$
We have an isomorphism
$$\langle - ,[S^{1}]\rangle:H^{2}_{Del}(S^{1},\Z)\stackrel{\cong}{\to} \C/\Z$$
given by evaluation.
  In \cite[(1.2),(1.3)]{MR1671481} Carey-Pincus  
 calculate the determinant invariant and joint torsion:
  $$d(\{u_{1},u_{2}\})=\tau(T_{u_{1}},T_{u_{2}})= \exp(2\pi i\langle u_{1}\cup u_{2},[S^{1}]\rangle)\ .$$   
 
 Combining this equality with \eqref{wefefwefewfewfewfewfewf} we get the equality
\begin{equation}\label{zgduiwegdwedwedwed89wed}\delta(b_{\Dirac}(\{u_{1},u_{2}\}))=\langle u_{1}\cup u_{2},[S^{1}]\rangle\ .\end{equation}
Note that we do not know whether the class $\{u_{1},u_{2}\}\in K_{2}(C^{\infty}(S^{1}))$ is topologically non-trivial.   
\bigskip

Using the multiplicative features of the differential regulator map $\hat \reg_{X}$ (see Remark \ref{ojefewlfjwelfjewfwef897}) one can also calculate $r_{\Dirac}(\sigma_{2}(\{u_{1},u_{2}\}))$ explicitly. The result is again the    right-hand side of \eqref{zgduiwegdwedwedwed89wed} as expected
by Conjecture \ref{ijdlkqwdwqdqwdqwd}. 
  We will not give the details  of the multiplicative theory  since it requires a set-up which is similar  to \cite{2013arXiv1311.1421B} but differs from the one used in the present paper.

 }
 \end{ex}

 \begin{rem}{\rm In this remark we recall the basic features of Deligne cohomology  used above. For $p\in \nat$   the Deligne cohomology group
$H^{p}_{Del}(X,\Z)$   is defined as the $p$'th hypercohomology of the complex of sheaves
$$0\to \underline{\Z}\to \Omega^{0}\to \Omega^{1}\to \dots\to \Omega^{p-1}\to 0\ .$$ We refer to \cite{MR827262} for a first definition and to \cite{MR1197353}, \cite[Sec. 3]{skript},  or \cite[4.3]{2013arXiv1311.3188B}) for  introductions to Deligne cohomology. 
In the original paper \cite{MR827262}  Deligne cohomology classes are called differential characters and a different grading convention was used.
We have a cup product
$$\cup:H^{p}_{Del}(X,\Z)\otimes H^{q}_{Del}(X,\Z)\to H^{p+q}_{Del}(X,\Z)$$ which turns Deligne cohomology into a graded commutative ring. Moreover, we have 
  a natural isomorphism of groups \begin{equation}\label{ghjghjdgwqduziqwdwqdwqdqwdwqd}H^{1}_{Del}(X,\Z)\cong C^{\infty}(X)^{\times}\ .\end{equation} Note that we get invertible complex-valued functions
  since 
   $\Omega^{*}$ is the  de  Rham complex of complex-valued forms.
  Finally, for a closed connected and oriented manifold $M$ of dimension $n-1$ we have an evaluation isomorphism
  $$\langle-, [M]\rangle :H^{n}_{Del}(M,\Z)\stackrel{\cong}{\to} \C/\Z\ .$$
   }\hB \end{rem}

 \begin{rem}\label{jwhdqkwdhqwd89qwd79w8q7d}{\rm
We refer to \cite[Appendix 4]{MR913964} for some information about the algebraic $K$-theory of the algebra of smooth functions on a manifold.

Let  $X$ be a compact manifold and $n\in \nat$ be odd. 
By  \cite[Thm A.4.6]{MR913964}  the rank of $$\im\left(K_{n}(C^{\infty}(X))\to K^{top}_{n}(C^{\infty}(X))\right)$$ is at least  $\dim H^{n}(X;\R)$. 
By  \cite[Thm A.4.6]{MR913964} this is true also if $X$ is oriented and $n=\dim(X)$ (not necessarily odd).
 
We have a decomposition 
$$K_{*}(C^{\infty}(X))\cong  K_{*}(\C)\oplus \tilde K_{*}(C^{\infty}(X))\ ,$$
where the first summand is induced by the inclusion $\C\to C^{\infty}(X)$ as constant functions, 
and the second summand 
is the kernel of the restriction to some point in $X$.
There are  non-trivial classes in $K_{d}(C^{\infty}(X))$ for arbitrary large odd $d\in \nat$. Note that classes coming from the summand $K_{d}(\C)$ are topologically trivial.

Theorem  \cite[Thm A.4.3]{MR913964} shows that the   group 
$\tilde K_{n}(C^{\infty}(X))$ itself is huge for $1\le n\le \dim(X)$.

We do not know whether there exists topologically non-trivial classes in degrees strictly larger than $\dim(X)$.

}\hB
\end{rem}

{\em Acknowledgement: The author  thanks J. Kaad for valuable discussions and clarification of the history of the  relation between joint torsion, determinant invariant, and the multiplicative character.}

\section{The construction of $\sigma_{d}$}\label{wdjqdjkwqdnkwqdwqd333}

\subsection{Elements of sheaf theory on manifolds}

Our main idea is to analyse the algebraic $K$-theory spectrum $\bK(C^{\infty}(X))$ of the algebra of complex-valued smooth functions on a manifold $X$
using the techniques of differential cohomology theory as developed in 
\cite{2013arXiv1311.3188B}. In the following we recall some of the basic notions.

Let $\Mf$ denote the site of smooth manifolds with corners with the open covering topology. 
 
 \begin{rem}{\rm An $n$-dimensional manifold with corners  is locally modeled by open subsets of  the subspace $[0,\infty)^{n}\subset\R^{n}$.
The category of manifolds with corners contains the unit  interval $[0,1]$, manifolds with boundary,
simplices $\Delta^{n}$. Furthermore, the category of manifolds with corners  is closed under taking products.

 }\hB \end{rem}
 
 \bigskip
 
 For a  presentable $\infty$-category $\cC$ (see \cite[Ch. 5]{HTT})
 we will consider the $\infty$-category of presheaves $\PSh_{\cC}(\Mf)$ and its full subcategory of sheaves $\Sh_{\cC}(\Mf)$ with values in $\cC$ on the site $\Mf$.
  \begin{ddd}A presheaf $G\in \PSh_{\cC}(\Mf)$ is a sheaf if for every manifold $M$ and every open covering
 $U\to M$ the natural map
\begin{equation}\label{hgh6769889898998}G(M)\to \lim_{\Delta} G(U^{\bullet})\end{equation} is
an equivalence.
\end{ddd} 
In this definition  the simplicial manifold $U^{\bullet}\in \Mf^{\Delta^{op}}$ is the \v{C}ech nerve of the open covering and the map \eqref{hgh6769889898998} is induced from the natural map $U^{\bullet}\to M$,
were $M$ is considered as a constant simplicial manifold.
By an application of the general theory \cite[6.2.2.7]{HTT} we get that  $\PSh_{\cC}(\Mf)$ and $\Sh_{\cC}(\Mf)$ are again presentable $\infty$-categories and that
there is an  adjunction 
$$L:\PSh_{\cC}(\Mf)\leftrightarrows \Sh_{\cC}(\Mf):\incl$$ 
between the inclusion of sheaves into presheaves and the sheafification functor $L$.

\bigskip

We use the unit interval  $I:=[0,1]$  in order to define the notion of homotopy invariance.
\begin{ddd} \label{kwhdwqjkdhkwdjwqddwqdwqdwqd32423423}A sheaf or presheaf $G$ on $\Mf$ is called homotopy invariant, if the map
$$G(M)\to G(I\times M)$$
induced by the projection $I\times M\to M$ is an equivalence for every smooth manifold $M$.
\end{ddd}
As in \cite[Sec. 2]{2013arXiv1311.3188B} one argues that
the full subcategories of homotopy invariant sheaves  $\Sh^{h}_{\cC}(\Mf)$ or homotopy invariant presheaves $\PSh^{h}_{\cC}(\Mf)$ are presentable. Their inclusions into  all sheaves or presheaves fit into adjunctions
\begin{equation}\label{dhqjdhkjqwdqwdw}
\cH^{pre}:\PSh_{\cC}(\Mf)\leftrightarrows \PSh^{h}_{\cC}(\Mf):\incl\ , \quad 
\cH:\Sh_{\cC}(\Mf)\leftrightarrows \Sh^{h}_{\cC}(\Mf):\incl\ .\end{equation}
The left adjoints are called homotopification functors. The homotopification functors for sheaves and presheaves  are related by the equivalence
\begin{equation}\label{xknksnxlsakxslkxasx879811}\cH\simeq L\circ \cH^{pre}\circ \incl\ ,\end{equation} see \cite[Prop. 2.6]{2013arXiv1311.3188B}.
 
\bigskip

Let $X$ be a smooth manifold. By
 $i_{X}:\Mf\to \Mf$ we denote the map of sites  given by $M\mapsto X\times M$.
It induces a pull-back
\begin{equation}\label{gdhdwhdkqwdqwdqwdwqdd}i^{*}_{X}:\PSh_{\cC}(\Mf)\to \PSh_{\cC}(\Mf)\ .\end{equation}
The functor $i_{X}^{*}$ has the following properties:
\begin{lem}\label{djlwqdqwdqwdwqdwd}{}
\begin{enumerate}
\item The functor $i_{X}^{*}$ preserves sheaves. \item The map $i_{X}^{*}$ preserves homotopy invariant presheaves and sheaves. 
\item For presheaves $i_{X}^{*}$ commutes with homotopification, i.e.  
the natural transformation   $\cH^{pre}\circ i_{X}^{*} \to i_{X}^{*}\circ \cH^{pre}$ is an equivalence.
\item If $X$ is compact, then the analogous statement holds true for sheaves, i.e. the natural map
$\cH\circ i_{X}^{*} \to i_{X}^{*}\circ \cH$ is an equivalence.
 
\end{enumerate}
\end{lem}
\proof
If $U\to M$ is an open covering covering of $M$, then $i_{X}(U)\to i_{X}(M)$ is  an open  covering of $X\times M$.   If $G$ is a presheaf, then the descent map of $i_{X}^{*}G$ with respect to 
$U\to M$ is the same as the descent map  of $G$ with respect to $i_{X}(U)\to i_{X}(M)$. This implies the first assertion.

A sheaf or presheaf $G$ is homotopy invariant by definition if the natural transformation $G\to i_{I}^{*}G$ (induced by the map $I\to *$) is an equivalence.
We have  equivalences of functors $i^{*}_{I}i^{*}_{X}\simeq i^{*}_{I\times X}\simeq i_{X}^{*}i_{I}^{*}$.
This implies that $i_{X}^{*}$ preserves homotopy invariant sheaves or presheaves. 

In order to see that
$i_{X}^{*}$ commutes with homotopification of presheaves we use the explicit formula for the homotopification given in \cite[Section 7]{2013arXiv1311.3188B}. We define the functor
 \begin{equation}\label{qhwdkqwhdkqwjdhwqdwqd78}\bs:\PSh_{\cC}(\Mf)\to \PSh_{\cC}(\Mf)\ , \quad \bs(G):=\colim_{\Delta^{op}} i_{\Delta^{\bullet}}^{*}G\ ,\end{equation} where $\Delta^{\bullet}$ is the cosimplicial manifold of standard simplices.
Then the homotopification on presheaves $\cH^{pre}$ is given by 
\begin{equation}\label{xknksnxlsakxslkxasx87981}\cH^{pre}\simeq \bs\ .\end{equation}   Since
the colimit for presheaves is taken object wise and $i_{\Delta^{\bullet}}^{*} i_{X}^{*}\simeq i_{X}^{*}i_{\Delta^{\bullet}}^{*}$ we see that the homotopification for presheaves commutes with $i_{X}^{*}$.

We now assume that $X$ is compact and $G$ is a sheaf. For $n\in \nat$ we let $S^{n}\in \Mf$ denote the $n$-dimensional sphere.
For every $n\in \nat$, using \cite[Prop. 7.6]{2013arXiv1311.3188B} at the marked places, we get the following chain of  equivalences: 
\begin{eqnarray*}\lefteqn{(i_{X}^{*}\cH G)(S^{n})\simeq(\cH G)(X\times S^{n})\stackrel{!}{\simeq} (\cH^{pre} G) (X\times S^{n})\simeq (i^{*}_{X} \cH^{pre} G)(S^{n})}&&\\&&\hspace{6cm}\simeq (\cH^{pre}i_{X}^{*}G)(S^{n})\stackrel{!}{\simeq} (\cH i_{X}^{*}G)(S^{n} )\ .\end{eqnarray*}
The Assertion 4. now follows from \cite[Lemma 7.3]{2013arXiv1311.3188B}
which states that an equivalence between objectes of $\Sh_{\cC}(\Mf)$ can be detected on the collection of spheres $S^{n}$, $n\in \nat$.
\hB  

\begin{rem}{\rm  Since $i_{X}^{*}$ preserves sheaves we have a natural transformation   \begin{equation}\label{gd8w7d987d987d98d}L\circ i_{X}^{*}\to i_{X}^{*}\circ L\ .\end{equation} In general it is not an equivalence.  For example, let $\cC:=\Ab$ and $\underline{\Z}^{pre}$ be the constant presheaf with value $\Z$ and $X$ consist of two points. 
Then we have
$i_{X}^{*}(L(\underline{\Z}^{pre}))(*)\cong \Z\oplus \Z$, but
$L(i_{X}^{*}\underline{\Z}^{pre})(*)\cong \Z$ and the map \eqref{gd8w7d987d987d98d} is the diagonal inclusion.
}\end{rem}

\bigskip

 In the present paper we will use the language of diffeological algebras. 
 \begin{ddd}\label{dkdljqwkldjwqldqwd}
A diffeological structure on an algebra $A$ over $\C$ is a 
subsheaf of algebras $A^{\infty}$   of the sheaf of algebras $M\mapsto \Hom_{\Set}(M,A)$
such that $A^{\infty}(*)=A$. A diffeological algebra is an algebra equipped with a diffeological structure. 
\end{ddd} 

\begin{rem}{\rm 
A sheaf $F$  of sets on $\Mf$ which is a  subsheaf of the sheaf of set-valued functions to $F(*)$ is also called a concrete sheaf. We refer to \cite{MR2805746} for a discussion of various variants of the definition of a diffeology. Our version is most similar to the notion of a Chen space, but not equal. A Chen space is a concrete sheaf of sets on the site of convex subsets with non-empty interior of euclidean spaces.    In contrast, our sheaves are defined on all manifolds with corners. 
}\hB
\end{rem}

\begin{ex}\label{dhqwkdhwqkdhqkjdqwdwq}{\rm
In the following we list some examples of diffeological algebras.
\begin{enumerate}
\item
The constant sheaf $\underline{A}$ generated by $A$ is the minimal diffeological structure, while the sheaf $M\mapsto \Hom_{\Set}(M,A)$ is the maximal diffeological structure on $A$.
\item
If $A$ is a diffeological algebra and $X$ is a smooth manifold, then we define the algebra
$C^{\infty}(X,A):=A^{\infty}(X)$. It has again a diffeological structure given by the sheaf
$i_{X}^{*}A^{\infty}$.
\item The algebra $\C$ has a diffeological structure such that $\C^{\infty}$ is the sheaf $M\mapsto C^{\infty}(M)$ of smooth $\C$-valued functions on $\Mf$. 
\item For a manifold $X$ we equip $C^{\infty}(X)$ with the diffeological structure defined in 2.
\item If $A$ is a locally convex algebra, then we have a natural notion of a smooth function $X\to A$.
The diffeological structure is given by $A^{\infty}(X):=C^{\infty}(X,A)$, where $C^{\infty}(X,A)$ denotes the algebra of smooth functions on $X$ with values in $A$. See Remark \ref{ldjlkqdjqwldqwdwqd7987} for more details.
\end{enumerate}
}\hB\end{ex}
\begin{rem}\label{ldjlkqdjqwldqwdwqd7987} {\rm 
In this remark we fix our conventions about smooth functions on manifolds with values in a locally convex algebra.
A locally convex vector space is a complex vector space whose topology is defined by a collection of seminorms. A locally convex algebra is a locally convex vector space such that the product induces a
continuous bilinear map
$A\times A\to A$.

A locally convex vector space has a natural uniform structure. Therefore the notions of completeness and completion are defined.

We now consider smooth functions with values in a locally convex vector space $A$ (see e.g. \cite[Sec. 40]{MR0225131}).
Let $U\subseteq \R^{n}$ open and consider a continuous function $f:U\to A$.
\begin{ddd} The function $f$  continuously differentiable  if there exists a continuous function $f^{\prime} :U\to   \Hom(\R^{n},A)$  such that for every seminorm $p$ on $A$ and every compact subset $K\subset U$ we have
$$\lim_{D\to 0} \sup_{x\in K}\:p\left(\frac{f(x+D)-f(x)-f^{\prime}(x)(D)}{\|D\|}\right)=0\ .$$
We call   $\partial_{i}f:=f^{\prime}(-)(e^{i})$ the partial derivatives of $f$ in the $i$'th direction.
We call $f$ smooth if it has all iterated continuous  partial derivatives.  
  \end{ddd}
We denote the iterated partial  derivatives by $f^{(k)}_{i_{1},\dots,i_{k}}$. We equip the complex vector space $C^{\infty}(U,A)$ with
the locally convex structure   determined by the seminorms
$$f\mapsto \sup_{x\in K} \: p( f^{(k)}_{i_{1},\dots,i_{k}}(x))\ .$$
The set of seminorms which generates the topology of $C^{\infty}(U,A)$ is thus indexed by compact subsets $K\subset U$, tuples $(i_{1},\dots,i_{k})$ of elements of $\{1,\dots,n\}$,  and  seminorms $p$ of $A$.  If $A$ is complete, then so is $C^{\infty}(U,A)$.

This definition of smooth $A$-valued functions extends to manifolds in a straightforward manner.

\bigskip

Let $X$ be a smooth manifold. Then the algebra $C^{\infty}(X,A)$ has two diffeological structures:
\begin{enumerate}\item
 The first comes from the construction 2. in Example \ref{dhqwkdhwqkdhqkjdqwdwq} above.
 \item 
 The second is induced from its locally convex structure. 
 \end{enumerate}
 These two structures coincide in view of the exponential law:
$$C^{\infty}(X,C^{\infty}(Y,A))\cong C^{\infty}(X\times Y,A)\ .$$
 
Since $C^{\infty}(M,A)$ is a subset of the set of all functions from  $M$ to $A$ it is clear that
these spaces of smooth functions for varying $M$ define a concrete sheaf, i.e. a diffeological structure on $A$. 

In the non-commutative geometry literature instead of $C^{\infty}(X,A)$ one often uses the projective tensor product   $C^{\infty}(X)\otimes_{\pi}A$. If $A$ is complete, then this gives an equivalent structure as we will explain below.
 We have a natural map
$$C^{\infty}(X)\otimes A\to C^{\infty}(X,A)$$ 
which is continuous with respect to projective topology on the algebraic tensor product.

In general, for locally convex vector spaces $V, W$ we let $V\otimes_{\pi}W$ denote the completion of the algebraic tensor product $V\otimes W$ with respect to the projective topology. If $A$ is a complete locally convex vector space,  then we get an isomorphism
$$C^{\infty}(X)\otimes_{\pi} A\stackrel{\cong}{\to} C^{\infty}(X,A)\ .$$
Here is a reference for this classical fact:
\begin{enumerate}
\item It follows from \cite[Thm 44.1]{MR0225131} that for a complete $A$  we have an isomorphism
$C^{\infty}(X)\otimes_{\epsilon}A\stackrel{\cong}{\to} C^{\infty}(X,A)$, where $\otimes_{\epsilon}$ denotes the completion of the algebraic tensor product in the $\epsilon$-topology.
\item It follows from  \cite[Thm 50.1]{MR0225131} that the locally convex vector space $C^{\infty}(X)$ is nuclear.
\item If one of the tensor factors is nuclear, then the natural map from the $\pi$- to the $\epsilon$-tensor product is an isomorphism by \cite[Thm 50.1]{MR0225131}. 
\end{enumerate}

 }\hB
\end{rem}

\begin{rem}{\rm In this remark we explain the relationship between the notions of homotopy invariance
according to Definition \ref{kwhdwqjkdhkwdjwqddwqdwqdwqd32423423} and diffeotopy invariance of functors defined on locally convex algebras
as considered e.g. in  \cite[Sec. 4.1]{MR2409415}.

Consider the category $\mathcal{L}oc\mathcal{A}lg_{1}$
of unital complete locally convex algebras.  We have a functor
$$\mathcal{L}oc\mathcal{A}lg_{1}\to \Sh_{\mathcal{L}oc\mathcal{A}lg_{1}}(\Mf)\ , \quad  A\mapsto A^{\infty}:=(M\mapsto C^{\infty}(M)\otimes_{\pi}A)\ .$$
Let $\cC$ be a presentable $\infty$-category.
 
\begin{lem} A functor $F: \mathcal{L}oc\mathcal{A}lg_{1}\to \cC$ 
is a diffeotopy invariant functor in the sense of \cite[Sec. 4.1]{MR2409415} if and only if
the presheaf $F(A^{\infty})\in \PSh_{\cC}(\Mf)$ is homotopy invariant in the sense of 
Definition \ref{kwhdwqjkdhkwdjwqddwqdwqdwqd32423423}.
\end{lem}
 \proof
 This is immediate from the definitions if one uses
the associativity of $\otimes_{\pi}$ and
$$C^{\infty}(I\times X)\cong C^{\infty}(I)\otimes_{\pi}C^{\infty}(X)\ .$$

}\hB\end{rem}

 \subsection{Algebraic $K$-theory and cyclic homology of smooth functions}
 
\subsubsection{Chain complexes and spectra} 
 
The purpose of this paragraph is to fix our conventions concerning chain complexes and spectra.
We fix some notation and introduce some basic constructions.
 
 Let $\Ch$ be the  category of (in general unbounded) chain complexes of abelian groups and chain morphisms. If $R$ is a ring, then we use $\Ch_{R}$ in order to denote the category of chain complexes of $R$-modules.

 We identify chain complexes with cochain complexes 
such that the chain complex
$$\dots \to C_{n+1}\to C_{n}\to C_{n-1}\to \dots$$
corresponds to the cochain complex
$$\dots \to C^{-n-1}\to C^{-n}\to C^{-n+1}\to \dots\ .$$

For an integer  $p\in \nat$ and a chain complex $(C,d)\in \Ch$   we  define  its shift by $p$ 
  by $C[p]^{n}:=C^{n+p}$. The differential of the shifted complex is given by   $(-1)^{p}d$.

 For $n\in \nat$ we let $H^{n}:\Ch\to  \Ab$ denote the $n$'th cohomology functor.
A chain map is a quasi-isomorphism if it induces an isomorphism in cohomology.
If we invert the quasi-isomorphisms in $\Ch$, then we get a stable $\infty$-category $\Ch[W^{-1}]$.
There are various ways to construct a model of  this $\infty$-category e.g. using model categories or
$dg$-enhancements. Since the constructions in the present paper are model independent
we will not discuss the details.

We have a natural functor $\iota:\Ch\to \Ch[W^{-1}]$. Our usual notation convention is  that  the italic letter $C$ denotes an object of $\Ch$, and the corresponding roman letter $\mathrm{C}$   denotes the object $\iota(C)\in \Ch[W^{-1}]$.
By the universal property of  $\Ch[W^{-1}]$ the cohomology functors descent to functors $H^{n}:\Ch[W^{-1}]\to \Ab$.

For an integer $p\in \Z$ and  a chain complex $C\in\Ch$
$$\dots\to C^{p-1}\to C^{p}\to C^{p+1}\to \dots$$    we  define its naive truncations $\sigma^{\ge p}C$ and $\sigma^{<p}C$   at $p$ by \begin{equation}\label{nnmasdhskdsad9898}\dots\to 0\to C^{p}\to C^{p+1}\to \dots\ , \quad \dots\to C^{p-2}\to C^{p-1}\to 0 \to\dots \ .\end{equation}
   We have 
    natural inclusion and projection morphisms \begin{equation}\label{ndjebdjkdjkqkdjwqdjhqwqwdwqdwq6d6}\sigma^{\ge p}C\to C\ , \quad C\to\sigma^{<p}C\ .\end{equation}

 
\begin{rem}{\rm  Note that  $\iota(\sigma^{\ge p }C)$ is well-defined, but  $\sigma^{\ge p}\iota(C)$ does not make sense. } \hB\end{rem}

By  $\Omega\in \Sh_{\Ch}(\Mf)$  we denote the sheaf of de Rham complexes on $\Mf$ of complex-valued differential forms. By \cite[Lemma 7.12]{2013arXiv1311.3188B}, 
for every $p\in \Z$   its truncation $\iota(\sigma^{\ge p}\Omega)$ is a sheaf, i.e 
 \begin{equation}\label{kksjwlkjslqwsqws}\iota(\sigma^{\ge p}\Omega)\in  \Sh_{\Ch[W^{-1}]}(\Mf)\ .\end{equation}   Note that $H^{p}(\iota(\sigma^{\ge p}\Omega))\cong \Omega^{p}_{cl}\in \Sh_{\Ab}(\Mf)$ is the sheaf of closed $p$-forms.

\bigskip

Let  $\Sp$ denote the stable  $\infty$-category of spectra. Again we will not discuss explicit models.
 For every $n\in \Z$ we have a functor
$\pi_{n}:\Sp\to \Ab$ which  maps a spectrum to its $n$'th homotopy group. The collection of these functors for all $n\in \Z$ detects equivalences in $\Sp$.

 We will frequently use the   Eilenberg-MacLane correspondence $H:\Ch[W^{-1}]\to \Sp$ 
 (see \cite[8.1.2.13]{HA}) which maps a chain complex to its associated Eilenberg-MacLane spectrum.  For $\mathrm{C}\in  \Ch[W^{-1}]$ we have the relations
$$\pi_{n}(H(\mathrm{C}))\cong H^{-n}(\mathrm{C})\ , \quad H(\mathrm{C}[p])\simeq \Sigma^{p} H(\mathrm{C})$$
between the homotopy groups of $H(\mathrm{C})$ and the cohomology groups of $\mathrm{C}$ on the one hand, and the
shifts by $p\in \Z$ in spectra and chain complexes, on the other.

 The Eilenberg-MacLane equivalence preserves limits. Hence it induces a map
$$H:\Sh_{\Ch[W^{-1}]}(\Mf)\to \Sh_{\Sp}(\Mf)$$ by objectwise application.
For example, by \eqref{kksjwlkjslqwsqws} we have the sheaf
\begin{equation}\label{kksjwlkjslqwsqws1}H(\iota(\sigma^{\ge p}\Omega))\in  \Sh_{\Sp}(\Mf)\ .\end{equation}
We  have  $\pi_{-p}(H(\iota(\sigma^{\ge p}\Omega)))\cong\Omega^{p}_{cl}$.

\subsubsection{Algebraic $K$-theory}

In this paragraph we call some basic facts from algebraic $K$-theory.
 We let  $\Alg$ denote  the category of associative unital algebras.   We have a functor $$\bK:\Alg\to \Sp$$ which maps an associative unital algebra to its connective algebraic  $K$-theory spectrum $\bK(A)$.   \begin{rem}{\rm
One way to construct this functor is as the following composition:
$$\bK(A):=\sp(\GrCompl(\Nerve(\Iso(\Proj(A)))))\ .$$
Here $\Proj(A)$ is the symmetric monoidal category of finitely generated projective $A$-modules with respect to the direct sum and $\Iso$ takes the underlying groupoid. The functor $\Nerve$ maps a symmetric monoidal category to its nerve which is  a commutative monoid in spaces. The group completion functor
$\GrCompl$ turns this monoid into a commutative group (i.e. an $E_{\infty}$-space) or equivalently into an infinite loop space.
Finally, the functor $\sp$ maps the infinite loop space to the corresponding spectrum.  We refer to \cite[Sec 6]{2013arXiv1311.3188B} and \cite[Sec. 6]{2013arXiv1306.0247B} for more details. In the present paper will not need any explicit construction of the algebraic $K$-theory functor.
}\hB \end{rem}

Let $\Alg_{\C} $ denote the  category of unital algebras over $\C$. We have a functor $$CC^{-}:\Alg_{\C}
\to \Ch$$ which maps an associative unital algebra $A$ over $\C$ to its negative
cyclic homology complex $CC^{-}(A)$.  
We define the negative cyclic homology of $A$ by 
\begin{equation}\label{jdjqhdqhdiuqwduwqidhqwdqwd89789737213213}HC^{-}_{*}(A):=H_{*}(CC^{-}(A))\ .\end{equation}
\begin{rem}\label{jefefelwfewfwf324324}{\rm 
For concreteness we will choose for $CC^{-}(A)$   the standard negative cyclic homology complex denoted by
$\mathrm{ToT} \:\cB C^{-}$ in \cite[5.1.7]{MR1600246}. Some constructions in the present paper will use this model explicitly.} \hB
\end{rem}
\item  We define the functor $\CC^{-}:\Alg_{\C}\to \Sp$ as the composition
$$\Alg_{\C}\stackrel{CC^{-}}{\to} \Ch\stackrel{\iota}{\to} \Ch[W^{-1}]\stackrel{H}{\to}\Sp\ .$$ 

We have a natural transformation of functors
\begin{equation}\label{djkejdkejdewd777wed}\ch^{GJ}:\bK\to \CC^{-}\ ,\end{equation} given by the Goodwillie-Jones Chern character.
 For the construction of the Goodwillie-Jones Chern character we refer to \cite{MR1275967} and
 \cite[Sec. 5]{MR2266961}. 

\subsubsection{Algebraic $K$-theory sheaves}

The goal of this paragraph is to introduce some basic notation which we will use throughout the rest of the paper. We consider  a diffeological algebra $A$ (see Definition \ref{dkdljqwkldjwqldqwd}) and form
the presheaf of spectra 
\begin{equation}\label{dhqwjdhwqjdhkqwdkd222}\check{\bK}_{A}\in \PSh_{\Sp}(\Mf)\ , \quad M\mapsto \bK(A^{\infty}(M))\ .\end{equation} 
 Its sheafification is a sheaf of spectra and will be denoted by 
\begin{equation}\label{deldqlkdwdwqddjh9879879}\bK_{A}:=L(\check{\bK}_{A})\in \Sh_{\Sp}(\Mf)\ .\end{equation}  

We apply these constructions to the diffeological algebra $\C$. We then
 have the equivalences
\begin{equation}\label{ghjwgdhjwgdqwjdqd666qwd}\check{\bK}_{C^{\infty}(X)}\simeq i_{X}^{*} \check{\bK}_{\C}\ ,  \quad \bK_{C^{\infty}(X)}\simeq L i_{X}^{*}\check{\bK}_{\C} 
\ .\end{equation}
The first follows from the definition of the diffeological structure on $C^{\infty}(X)$, and the second is then  a reformulation of the definition above.

\bigskip

Applying the negative cyclic homology complex to the sheaf $\C^{\infty}$
 we obtain the presheaf of chain complexes 
$$CC^{-}(\C^{\infty})\in \PSh_{\Ch}(\Mf)\ , \quad M\mapsto CC^{-}(C^{\infty}(M))\ .$$
\begin{rem}{\rm
Note that we do not complete or sheafify the tensor products involved in the definition of the negative cyclic homology complex, but see Remark \ref{djwd444jwqdwqdwq}.  
}\hB \end{rem}

In the following we define a differential geometric analog $DD^{-}$  of $CC^{-}(\C^{\infty})$ and a comparison map
$$\pi^{-}:CC^{-}(\C^{\infty})\to DD^{-}\ .$$ 
%

\begin{ddd}\label{dhwqdhqwkdqwdqwd}We define the sheaf of chain complexes $DD^{-}\in \Sh_{\Ch}(\Mf)$
$$ DD^{-}:= \prod_{p\in \Z} DD^{-}(p)\ , \quad DD^{-}(p):=  (\sigma^{\ge p} \Omega)[2p]\ .$$
 We further define a map of presheaves of chain complexes
\begin{equation}\label{kkckaslckjaslkcjaslkcjasc}\pi^{-}:CC^{-}(\C^{\infty})\to DD^{-}\end{equation} by 
\begin{eqnarray*}\lefteqn{CC^{-}(\C^{\infty}(X))_{q,p}\ni f_{0}\otimes \dots \otimes f_{p-q}\mapsto \frac{1}{(p-q)!}f_{0}df_{1}\wedge \dots df_{p-q}}&&\\&&\hspace{8cm}\in F^{p-q}\Omega^{p-q}(X) \subset     DD^{-}(p)(X)^{-p-q} \ .\end{eqnarray*}
\end{ddd}

Here the index $(\dots)_{q,p}$ refers to the component of the bicomplex $\cB C^{-}$ in \cite[5.1.7]{MR1600246}.
Using the formulas  \cite[Sec.2.3.2]{MR1600246} we conclude that $\pi^{-}$ 
is compatible with the differentials.
\bigskip
\begin{rem}\label{djwd444jwqdwqdwq}{\rm 
For a manifold $X$ let $CC^{cont,-}(C^{\infty}(X))$ be the analog of $CC^{-}(\C^{\infty}(X))$ defined using completed (but not sheafified) tensor products. Then we have a factorization of $\pi^{-}$:
\begin{equation}\label{}  CC^{-}(\C^{\infty}(X))\to CC^{cont,-}(C^{\infty}(X))\stackrel{\pi^{cont,-}}{\to} DD^{-}(X)\ .\end{equation} 
 The second map is  quasi-isomorphism by the well-known calculation of the continuous negative cyclic homology of the algebra of smooth functions on a smooth manifold.
 We will use the continuous version of cyclic homology and $\pi^{cont,-}$  in Subsection \ref{dwlkdwqldqwdqwdqwdwdwd} below.
 }\hB \end{rem}
 


%
We further define the presheaves  
 \begin{equation}\label{sqwmsqnwsswsqwssqw798}\CC^{-}_{\C}:=  H\circ \iota\circ  CC^{-}(\C^{\infty}) \ , \quad \mathrm{DD}^{-}:=\iota(DD^{-})\ ,\quad \DD^{-}:= H \circ \mathrm{ DD}^{-}\ .\end{equation}
The $\Ch[W^{-1}]$-valued presheaf  $\mathrm{DD}^{-}$ is a sheaf   by \eqref{kksjwlkjslqwsqws1}.
As remarked above, the Eilenberg-MacLane functor $H$ preserves sheaves. Therefore $\DD^{-}$ is a $\Sp$-valued sheaf.

\bigskip

By its naturality the  Goodwillie-Jones Chern character \eqref{djkejdkejdewd777wed} provides a map 
$$\ch^{GJ}:\check{\bK}_{\C}\to \CC^{-}_{\C}\ .$$
between presheaves of spectra.   
\begin{ddd}
We define the regulator morphism $\check{\reg}$ of presheaves of spectra as the composition
 \begin{equation}\label{edededdduuizwdqwd}\check{\reg}:\check{\bK}_{\C} \stackrel{\ch^{GJ}}{\to}  \CC^{-}_{\C}    \stackrel{\pi^{-}}{\to} \DD^{-} \ .\end{equation}
Furthermore, for a smooth manifold $X$,  
 we define  the morphism of sheaves of spectra
 \begin{equation}\label{dhqwhdqwdwqd}\reg_{X}:\bK_{C^{\infty}(X)} \stackrel{\eqref{ghjwgdhjwgdqwjdqd666qwd}}{\simeq} Li_{X}^{*}\check{\bK}_{\C}\stackrel{Li_{X}^{*}\check{\reg}}{\to} Li_{X}^{*}\DD^{-} \simeq i_{X}^{*}\DD^{-}\end{equation}
\end{ddd} 
The last equivalence in \eqref{dhqwhdqwdwqd} follows from the fact that $\DD^{-}$ and therefore $i_{X}^{*}\DD^{-}$ are  sheaves. \hB
 
\begin{rem}
{\rm In this paper we usually call transformations from $K$-theory to cyclic homology Chern characters, and transformations from  $K$-theory  to differential forms regulators. There is one exception, namely the usual Chern character from topological $K$-theory to cohomology with complex coefficients calculated by the de Rham cohomology.

}
\end{rem}

\subsection{Homotopification and regulator maps}\label{jckjcnkwcewcec}

For a sheaf   $G$  with values in a stable $\infty$-category  (e.g. $\Ch[W^{-1}]$ or $\Sp$) we have a functorial homotopification fibre sequence of sheaves  
\begin{equation}\label{jjkfewfwefewfewffewf4333}\cA(G) \to  G\to \cH(G)\to \Sigma \cA
(G)\ ,\end{equation}
see \cite[Def. 3.1]{2013arXiv1311.3188B}.  The map  $G\to \cH(G)$ is the unit of the homotopification functor $\cH$ introduced in \eqref{dhqjdhkjqwdqwdw}, and $\cA$ by definition takes the fibre of this unit. The sheaf $G$ is homotopy invariant if and only if $\cA(G)\cong 0$. 

\bigskip

Let $A$ be a diffeological algebra (Definition \ref{dkdljqwkldjwqldqwd}) and $\bK_{A}$ be as in \eqref{deldqlkdwdwqddjh9879879}.
\begin{ddd}\label{dhqkdhwqkjdqwdqwdqwdwqd23424}
We define the sheaves of spectra $$\bK_{A}^{top}:=\cH(\bK_{A})\ , \quad \bK^{rel}_{A}:=\cA(\bK_{A})\ .$$
We call the evaluations $\bK^{top}(A):=\bK_{A}^{top}(*)$ and $\bK^{rel}(A):=\bK_{A}^{rel}(*)$
the topological and relative $K$-theory spectra of $A$.
\end{ddd}
 Note that the topological and relative  $K$-theory spectra depend on the diffeological structure on $A$. They fit into the fibre sequence of spectra
 \begin{equation}\label{lkdjlkedjwedwed}\bK^{rel}(A)\to \bK(A)\to \bK^{top}(A)\to \Sigma \bK^{rel}(A)\end{equation}
derived from \eqref{jjkfewfwefewfewffewf4333} by evaluation at $*$.

\begin{rem}\label{kdjqkldjqlwkdjwqdwqdqwdwqdwd}{\rm

A Fr\'echet algebra has a natural diffeological structure such that $A^{\infty}(M)$ is the algebra of smooth maps $M\to A$. In this case our definition
of $\bK^{top}(A)$ coincides with that  given in \cite[Sec. 3.1]{MR972606}.
Indeed, in this reference the authors apply Quillen's $+$-construction to the classifying space 
of the simplicial group $GL(C^{\infty}(\Delta^{\bullet},A))$. Using \eqref{qhwdkqwhdkqwjdhwqdwqd78} we can identify the resulting space with $\Omega^{\infty}(\bs(\bK_{A})(*))$. 
We now use the equivalence $\bs(\bK_{A})(*)\simeq \bK_{A}^{top}(*)$ which follows from the combination of
\eqref{xknksnxlsakxslkxasx87981} and \eqref{xknksnxlsakxslkxasx879811}. 
As a consequence,  the relative $K$-theory of a
Fr\'echet algebra  defined in \cite[Sec. 3.1]{MR972606} is isomorphic to our version.

More generally, if $A$ is a complete locally convex algebra, then in \cite[Def. 4.1.3]{MR2409415}
the notion of the diffeotopy $K$-theory spectrum was defined. This definition is just   
$\bs(\bK_{A})(*)$ written down in different symbols. Therefore our topological or relative  $K$-theory of a complete locally convex algebra also coincides with the diffeotopy or relative $K$-theory of 
 \cite{MR2409415}.


} \hB\end{rem}

The homotopification of the sheaf $\mathrm{DD}^{-}\in \Sh_{\Ch[W^{-1}]}(\Mf)$ defined in \eqref{sqwmsqnwsswsqwssqw798}
can be calculated explicitly  again in terms of differential forms. To this end we  introduce  the two-periodic de Rham complex.
 \begin{equation}\label{kkjlkwqjslqwjswsqwsw}DD^{per}:=\prod_{p\in \Z} \Omega[2p]\in \Sh_{\Ch}(\Mf)\ .\end{equation}  Its cohomology  \begin{equation}\label{kqwjdlkqwdjklqwjdqwdqwdq}HP^{*}(X):=H^{*}(DD^{per}(X))\end{equation}
is called  the periodic cohomology of the manifold $X$. The periodicity is implemented by the shift isomorphism, which for 
 $k\in \Z$  is given by 
\begin{equation}\label{shiftiso}\iota_{2k}:HP^{*}(X)\stackrel{\cong}{\to}  HP^{*+2k}(X)\ , \quad \iota_{2k}([\omega(p)])_{p\in \Z}=([\omega(p-k)])_{p\in \Z}\ .\end{equation}
We further
  define  
\begin{equation}\label{hdgjhqwgdwqgdhjwqdghjwdwd234324324}\mathrm{DD}^{per}:=\iota(DD^{per})\in \Sh^{h}_{\Ch[W^{-1}]}(\Mf)\ .\end{equation}
A priori we have $\mathrm{DD}^{per}\in \PSh_{\Ch[W^{-1}]}(\Mf)$.
In order to see that $\mathrm{DD}^{per}$ is a sheaf we use \cite[Lemma 7.12]{2013arXiv1311.3188B}. Moreover, 
since $ \Omega $ resolves the constant sheaf $\underline{\C}$, the sheaf $\iota(\Omega)$ is homotopy invariant. Consequently, the sheaf $\mathrm{DD}^{per}$ is homotopy invariant, too.

\bigskip 
 In view of Definition \ref{dhwqdhqwkdqwdqwd} of $DD^{-}$ and \eqref{ndjebdjkdjkqkdjwqdjhqwqwdwqdwq6d6}  we have a natural inclusion of sheaves of chain complexes $DD^{-}\to DD^{per}$.
\begin{lem}\label{khdlqwdjlqwjdqwdqwdqwdqwdwqd}
The induced map $ \mathrm{DD}^{-}\to \mathrm{DD}^{per}$ is equivalent to the homopification morphisms of $\mathrm{DD}^{-}$. In particular we have an equivalence $\cH(\mathrm{DD}^{-})\simeq \mathrm{DD}^{per}$. 
\end{lem}
\proof
By \cite[Lemma 7.15]{2013arXiv1311.3188B} we know that the inclusion
$\iota(\sigma^{\ge p}\Omega) \to \iota(\Omega)$ is equivalent to the homotopification map. This implies  that the natural inclusion $\mathrm{DD}^{-}\to \mathrm{DD}^{per}$ is equivalent to the homotopification map $\mathrm{DD}^{-}\to \cH(\mathrm{DD}^{-})$. \hB

We now provide  chain complex model for $\cA(\mathrm{DD}^{-})$.
We  define the sheaf of chain complexes
\begin{equation}\label{mdmwwqdqwd87qw9}DD:=\prod_{p\in \Z} DD(p) \in \Sh_{\Ch}(\Mf)\ , \quad DD(p):=(\sigma^{\le p}\Omega)[2p] \ .\end{equation}
It fits into  exact sequence of sheaves of chain complexes
\begin{equation}\label{qdhwgdqwd6w78d6z}0\to DD^{-}\to DD^{per}\to DD[2]\to 0\ .\end{equation}
The second map in this sequence   is induced by the family of maps of chain complexes
\begin{equation}\label{jkldjwqdlkwqjdqwdwdqwd}DD^{per}(p)\cong \Omega[2p]\stackrel{\eqref{ndjebdjkdjkqkdjwqdjhqwqwdwqdwq6d6}}{\to} (\sigma^{<p}\Omega)[2p]\cong ((\sigma^{\le p-1}\Omega)[2(p-1)])[2]\cong DD(p-1)[2]\ .\end{equation}
By  \cite[Lemma 7.12]{2013arXiv1311.3188B} the object
$\mathrm{DD}:=\iota(DD)$ is a sheaf with values in $\Ch[W^{-1}]$. 
From \eqref{qdhwgdqwd6w78d6z} we get a fibre sequence
\begin{equation}\label{klewffwefewfewfewfewf}\dots\to \mathrm{DD}[1]\to \mathrm{DD}^{-}\to \mathrm{DD}^{per}\to  \mathrm{DD}[2]\to \dots\ .\end{equation}
Lemma \ref{khdlqwdjlqwjdqwdqwdqwdqwdwqd} implies that this sequence is equivalent to the homotopification sequence  \eqref{jjkfewfwefewfewffewf4333} applied to $\mathrm{DD}^{-}$.
\begin{kor}
We have  an equivalence $  \cA(\mathrm{DD}^{-})\simeq \mathrm{DD}[1] $.\end{kor}
We define the sheaves of spectra $$\DD^{per}:=H(\mathrm{DD}^{per})\ , \quad \DD:=H(\mathrm{DD}) \ .$$ 
If we apply $H $ to the sequence \eqref{qdhwgdqwd6w78d6z}, then we get the fibre sequence of spectra
\begin{equation}\label{kdkjdkqwjdklwjqdlkqwdqwdwqd}\Sigma \DD\to \DD^{-}  \to  \DD^{per}\to \Sigma^{2}\DD\ .\end{equation}
\begin{kor}
The fibre sequence of spectra  \eqref{kdkjdkqwjdklwjqdlkqwdqwdwqd}
  is equivalent to the homotopification sequence \eqref{jjkfewfwefewfewffewf4333} applied to  $\DD^{-}$. \end{kor}

\bigskip

Let now $X$ be a manifold.
If we apply the homotopification sequence \eqref{jjkfewfwefewfewffewf4333}  to the morphism $\reg_{X}:\bK_{C^{\infty}}\to i_{X}^{*}\DD^{-}$, then we get the first two lines of the following diagram: \begin{equation}\label{23ek23e3e23e}
\xymatrix{\bK_{C^{\infty}(X)}^{rel}\ar@/^1.2cm/@{..>}[dd]^(.3){\reg_{X}^{rel}}\ar[d]\ar[r]^{\partial}&\bK_{C^{\infty}(X)}\ar[d]^{\reg_{X}} \ar[r]&\bK_{C^{\infty}(X)}^{top}\ar[d]\ar@/^1.2cm/@{..>}[dd]^(.3){\reg_{X}^{top}}\ar[r]&\ar[d]\Sigma\bK_{C^{\infty}(X)}^{rel}\ar@/^1.2cm/@{..>}[dd]^(.3){\reg_{X}^{rel}}\\\cA(i_{X}^{*}\DD^{-}) \ar@{-->}[d]\ar[r]&i_{X}^{*}\DD^{-}\ar[r]\ar@{=}[d]&\cH(i_{X}^{*}\DD^{-})\ar[r]\ar@{-->}[d]^{1}&\Sigma\cA(i_{X}^{*}\DD^{-}) \ar@{-->}[d]\\ i_{X}^{*}\Sigma \DD\ar[r]&i_{X}^{*}\DD^{-}\ar[r]&i_{X}^{*}\DD^{per}\ar[r]&i_{X}^{*}\Sigma^{2}\DD} \ .
\end{equation} 
Since $i_{X}^{*}\DD^{per}$ is homotopy invariant we obtain the dashed morphism denoted by $1$ and the lower middle square from the universal property of the homotopification. The lower part of the diagram is now defined 
by extension of the lower middle square to a morphism of fibre sequences.
We use the diagram in order to define
the relative and and topological regulator maps $\reg^{rel}_{X}$ and $\reg^{top}_{X}$ as indicated.

\begin{rem}\label{kjdqlwkdlqwdqwldlqwd}{\rm
 If $X$ is a compact manifold, then by Lemma \ref{djlwqdqwdqwdwqdwd}. 4. we have   equivalences
 \begin{equation}\label{dlqwjdlwqdjqwoldwqoldkwqd676}i_{X}^{*}\circ \cA\simeq \cA\circ i_{X}^{*}\ , \quad i_{X}^{*}\circ \cH\simeq \cH\circ i_{X}^{*}\ .\end{equation}
If we apply $i_{X}^{*}$ to \eqref{kdkjdkqwjdklwjqdlkqwdqwdwqd}, then the resulting sequence    \begin{equation}\label{kdkjdkqwjdklwjqdlkqwdqwdwqd111}i_{X}^{*}\Sigma \DD\to i_{X}^{*} \DD^{-}  \to  i_{X}^{*} \DD^{per}\to  i_{X}^{*}\Sigma^{2}\DD\end{equation} is equivalent to the homotopification sequence of   $i_{X}^{*}\DD^{-}$.
In particular, for a compact manifold $X$ the lower two lines in \eqref{23ek23e3e23e} are equivalent, and we have the equivalences
$\reg_{X}^{rel}\simeq \cA(\reg_{X})$ and $\reg_{X}^{top}\simeq \cH(\reg_{X})$.
}\hB\end{rem}



By its definition \eqref{dhqwhdqwdwqd} we have the following factorization of the regulator $\reg_{X}$:
 \begin{equation}\label{ddewdewdewldkllklk}
   \bK_{C^{\infty}(X)} \stackrel{\mbox{\tiny def}}{=} L i_{X}^{*} \check{\bK}_{\C} \stackrel{\eqref{gd8w7d987d987d98d}}{\to} i_{X}^{*}L\check{\bK}_{\C}  \stackrel{\mbox{\tiny def}}{=} i_{X}^{*}\bK_{\C} \to  i_{X}^{*} \DD^{-}\ . 
\end{equation} Using the universal property of the homotopification morphism
$\bK_{C^{\infty}(X)}\to \bK^{top}_{C^{\infty}(X)}$ we get the factorization of $\reg_{X}^{top}$: \begin{equation}\label{jkwhdkqwdqwdwqd}
\xymatrix{\bK^{top}_{C^{\infty}(X)}\ar[r]^{!!!}\ar@/^1cm/@{..>}[rr]^{\reg^{top}_{X}}&i_{X}^{*} \bK_{\C}^{top}\ar[d]^{\simeq}\ar[r]^{i_{X}^{*} \reg^{top}}&i_{X}^{*}\DD^{per}\\&i_{X}^{*}\underline{\ku}\ar[ur]_{i_{X}^{*}\ch^{gj}}&}\ .
\end{equation}

The  identification of $\bK^{top}_{\C}\simeq \underline{\ku}$ follows from \cite[Lemma 6.3]{2013arXiv1311.3188B}, where $\underline{\ku}$ is the constant sheaf of spectra generated by  the connective topological $K$-theory spectrum  
  $\ku\simeq \bK^{top}_{\C}(*)$ of $\C$.
We define $ \ch^{gj}:\underline{\ku}\to \DD^{per}$ so that the diagram commutes.

 \bigskip

If we evaluate $i_{X}^{*}\reg^{top}$ at a point, identify $(i_{X}^{*}\bK^{top}_{\C})(*)\simeq \underline{\ku}(X)$, and take homotopy groups, then 
we get a homomorphism (natural in $X$)
$$\ch^{gj}:\mathbf{ku}^{*}(X)\to \pi_{-*}(\DD^{per}(X))\stackrel{\eqref{kqwjdlkqwdjklqwjdqwdqwdq}}{\cong} HP^{*}(X)\ .$$
  The origin of this map is the Goodwillie-Jones Chern character $\ch^{GJ}$, see \eqref{djkejdkejdewd777wed}.
In this situation we also have the usual Chern character defined by Chern-Weil theory
$$\ch^{cw}:\mathbf{ku}^{*}(X)\to  HP^{*}(X)\ .$$
We refer to \cite[Sec. 6.1]{2013arXiv1311.3188B} for a construction of the Chern character $\ch^{cw}$  using differential cohomology methods.
The following lemma is probably well-known. For completeness we sketch an argument.
\begin{lem}\label{jjljldkqjwdqwdqwd}
We have  the equality of Chern character maps $\ch^{gj}=\ch^{cw}:\ku^{*}(X)\to HP^{*}(X)$. 
\end{lem}
\proof
Since both Chern characters arise from maps between homotopy invariant sheaves of spectra
they are characterized by their evaluation at the point. Since the target is rational
they are determined by their actions on homotopy groups. Hence it suffices to check
that
$$\ch^{gj}=\ch^{cw}:\ku^{*}(*)\to HP^{*}(*)\ .$$
One can now use the fact that both Chern characters are multiplicative (see \cite{MR1275967} for the multiplicitivity of $\ch^{GJ}$) in order to reduce
to the cases $*=0$ and $2$.
For $*=0$ one checks the equality directly going through the definitions.
For $*=2$ we argue as follows. We know by an explicit calculation that
$$\ch^{GJ}:\bK_{1}(C^{\infty}(S^{1}))\to H_{1}(CC^{cont,-}(C^{\infty}(S^{1})))\cong \Omega^{1}(S^{1})$$ maps the class $\iota(u)\in \bK_{1}(C^{\infty}(S^{1}))$
  of a unit $u\in C^{\infty}(S^{1})^{\times}$  to the form $\frac{1}{2\pi i}d\log u\in \Omega^{1}(S^{1})$.

We now consider specifically the embedding $u:S^{1}\to \C^{*}$.
The class $\iota(u)$ is then mapped to the generator  $x\in \ku^{-1}(S^{1})\cong \Z$
under the composition
$$ \bK_{1}(C^{\infty}(S^{1}))\to  \bK^{top}_{1}(C^{\infty}(S^{1}))\xrightarrow{\eqref{jkwhdkqwdqwdwqd},!!!} \pi_{1}(\bK^{top}_{\C}(S^{1}))\cong  \ku^{-1}(S^{1})\ .$$
Therefore $\ch^{gj}(x)\in HP^{-1}(S^{1})$ is given by the class
$(c(p))_{p\in \Z}\in H^{-1}(DD^{per}(S^{1}))$ with 
$$c(p):=\left\{\begin{array}{cc} [\vol_{S^{1}}]\in H^{-1}(DD^{per}(p)(S^{1}))   &p=1\\[0.3cm]0&else\end{array}\right.\ ,$$
where $\vol_{S^{1}}=\frac{1}{2\pi i} d\log u$ is the normalized volume form of $S^{1}$.
It follows by suspension that $\ch^{gj}: \mathbf{ku}^{-2}(*) \to  HP^{-2}(*)$ maps the generator
of $\ku^{-2}(*)\cong \Z$ to the class 
 $(b(p))_{p\in \Z}\in H^{-2}(DD^{per}(*))$
given by $$b(p):=\left\{\begin{array}{cc}[1]\in H^{-2}(DD^{per}(p)(*))   &p=1\\[0.3cm]0&else\end{array}\right.\ .$$
The same holds true for $\ch^{cw}$.  
\hB

\subsection{The differential regulator map}

In this subsection  we introduce a version $\hat \ku$ of the Hopkins-Singer  differential cohomology   associated to $\mathbf{ku}$ and a differential regulator  map $$\hat \reg_{X}:\bK(C^{\infty}(X))\to \hat{ \mathbf{ku}}(X)\ .$$
\bigskip

 \begin{ddd}\label{fljfkelwjflewfewfef}
We define the  Hopkins-Singer differential connective complex  $K$-theory $ \hat{\mathbf{ku}} \in \Sh_{\Sp}(\Mf)$ by
 \begin{equation}\label{sjaskjd8jkhdkqwdqw} \xymatrix{\hat{\mathbf{ku}} \ar[r]^{R}\ar[d]^{I}& \DD^{-}\ar[d]\\\underline{\mathbf{ku}} \ar[r]^{\ch^{cw}}&   \DD^{per}}\ .\end{equation} We further define the  differential  connective complex   $K$-theory groups of a manifold  $X$ by
 $$\widehat{ku}^{-*}(X):=\pi_{*}(\hat\ku (X))\ .$$
 \end{ddd}
 For a general discussion of Hopkins-Singer differential cohomology theories we refer to \cite[Sec. 4.4]{2013arXiv1311.3188B}. The map   $I$ takes the underlying  $\ku$-class, and the map   $R$ is called the curvature morphism.
 
 \bigskip
 
We fix a manifold $X$ and    recall the definition \eqref{dhqwhdqwdwqd} of $\reg_{X}$ and  the (Chern-Weil version) of the  Chern character $\ch^{cw}$   which is equivalent to $\ch^{gj}$ by Lemma \ref{jjljldkqjwdqwdqwd}.

 The   middle square of the diagram \eqref{23ek23e3e23e}  together with \eqref{ddewdewdewldkllklk} and  Lemma \ref{jjljldkqjwdqwdqwd}  gives a  commutative square 
   \begin{equation}\label{kwndkqwldwqdwqdd}\xymatrix{\bK_{C^{\infty}(X)}\ar[d]\ar[r]^{\reg_{X}}&\ar[d] i_{X}^{*}\DD^{-}  \\
 i_{X}^{*}\underline{\mathbf{ku}}\ar[r]^{i_{X}^{*}\ch^{cw}}&i_{X}^{*}\DD^{per}& \ .}  \end{equation}
 
 \begin{ddd}\label{dkjqldkjqwldqwdwqdw4} For a smooth manifold $X$
 we define the differential regulator map $$\hat \reg_{X}:\bK_{C^{\infty}(X)}\to i_{X}^{*}\hat{\mathbf{ku}}$$
 as the map between sheaves of spectra naturally induced by the square \eqref{kwndkqwldwqdwqdd} and the universal property of the pull-back square \eqref{sjaskjd8jkhdkqwdqw}.
 The evaluation of $\hat \reg_{X}$ at  a point  gives a map of spectra
 $$\hat \reg_{X}:  \bK (C^{\infty}(X))\to \hat{\mathbf{ku}} (X)\ .$$
   \end{ddd}

\begin{rem}\label{ojefewlfjwelfjewfwef897}
{\rm For commutative algebras it is possible to refine the Goodwillie-Jones Chern character $\ch^{GJ}$ to a morphism between commutative ring spectra. The spectra occuring on the right   corners of the diagrams  \eqref{sjaskjd8jkhdkqwdqw} and \eqref{kwndkqwldwqdwqdd} are obtained by an application of $H\circ \iota$ to sheaves of commutative differential graded algebras and therefore are commutative ring spectra, as well. Since the morphisms $\reg_{X}$ and $\ch^{cw}$ and the diagram \eqref{kwndkqwldwqdwqdd} have refinements to morphisms between sheaves of commutative ring spectra, 
the differential regulator   naturally becomes
a morphism between commutative ring spectra, too.

\bigskip

The multiplicative structure is helpful if one wants to calculate $\sigma_{d}$
of products like $$\iota(u_{1})\cup \dots\cup \iota(u_{d})\in K_{d}(C^{\infty}(X))$$ for a collection of invertible functions $(u_{i})_{i=1,\dots,d}$ in $C^{\infty}(X)$. Since our main results do not use the multiplicative structure, in the present paper we will not discuss its details further.
} \hB 
\end{rem}

\begin{rem}{\rm In general,
for an algebraic   $K$-theory class $x\in K_{d}(C^{\infty}(X))$ the differential $K$-theory class $\hat \reg_{X}(x)\in \hat \ku^{-d}(X)$ contains  strictly more information  than just the underlying topological class $ x^{top} :=I(\hat \reg_{X}(x))\in \ku^{-d}(X)$ and
the regulator $\reg_{X}(x)=R(\hat \reg_{X}(x))\in H^{-d}(DD^{-}(X))$. It is essentially this  additional secondary information which we   detect by the map $\sigma_{d}$ in Theorem \ref{djekldjlwedwedewdew44dewde}.

}
\end{rem}

 \subsection{The construction of $\sigma_{d}$}

We consider an integer $d\in \nat$ and a   smooth manifold $X$ such that its dimension satisfies $\dim(X)\le d-1$.   The goal of the present subsection is to construct the homomorphism
$$\sigma_{d}:K_{d}(C^{\infty}(X))\to \ku\C/\Z^{-d-1} (X)$$ asserted in Theorem \ref{djekldjlwedwedewdew44dewde}.
The idea is to obtain this map by specializing the differential regulator map $\hat \reg_{X}$.
%

\bigskip

For the moment let $X$ be an arbitrary smooth manifold.   
 We fix an integer $d\in \Z$. In the following definition the symbol
$R$ denotes the curvature map  the differential cohomology theory $\hat \ku$, see   \eqref{sjaskjd8jkhdkqwdqw}.
\begin{ddd}We define the subgroup of flat classes
$$\widehat{ku}^{-d}_{flat}(X):=\ker\left(R:\widehat{ku}^{-d}(X)\to \pi_{d}(\DD^{-}(X)) \right)\ .$$
 \end{ddd}

The flat subgroup of a Hopkins-Singer differential cohomology theory can be calculated explicitly.

\begin{lem}\label{kjewklfjwelkfewfwef}
 If $\dim(X)\le d$, then
we  have a natural isomorphism
\begin{equation}\label{dwekldjewlkdjweldjewdewd7we8d7wed} \widehat{ku}^{ -d}_{flat}(X)\cong \ku\C/\Z^{-d-1}(X) \ .\end{equation}
 \end{lem}
\proof
  The complex of sheaves of abelian groups $\Omega$ is a resolution of the constant sheaf $\underline{\C}$.  Consequently $DD^{per}$ resolves the constant sheaf $ \underline{\prod_{p\in \Z}\C[2p]}$. We thus obtain an equivalence
\begin{equation}\label{jskljkjlkjwsqws67567572613}\DD^{per}\simeq  \underline{H(\iota(\prod_{p\in \Z}\C[2p]))}\ .\end{equation}

 We have a fibre sequence of pull-back squares
$$\dots \to\xymatrix{\Sigma^{-1}\DD^{-}\ar@{=}[r] \ar[d]&\Sigma^{-1}\DD^{-}\ar[d]\\ 0\ar[r]&0}\to   \xymatrix{   \underline{\mathbf{F}} \ar[r]  \ar[d]&0\ar[d]\\\underline{\ku}\ar[r]&\DD^{per} }\to \xymatrix{\hat \ku\ar[r] \ar[d]&\DD^{-}\ar[d]\\\underline{\ku}\ar[r]&\DD^{per}}\to \dots\ ,$$
 
We calculate the fibre $\mathbf{F}$ of \begin{equation}\label{kwdnwqlkdqwdqwd}\ch^{cw}:\ku\to  H(\iota(\prod_{p\in \Z}\C[2p]))\ .\end{equation}   
%
%
Using the decomposition
$$\prod_{p\in \Z}\C[2p]\cong \prod_{p\in \nat}\C[2p]\oplus\prod_{p=1}^{\infty} \C[-2p] $$
and the fact that $\ku$ is connective
we get the decomposition of $\mathbf{F}$ into the sum of the fibres of the morphisms of spectra
$$\ku\to H(\iota(\prod_{p\in \nat}\C[2p]))\ , \quad 0\to H(\iota(\prod_{p=1}^{\infty} \C[-2p]))\ .$$
The fibre of the first morphism is equivalent to  $\Sigma^{-1} \ku\C/\Z$, and the fibre of the second is
$\Sigma^{-1} H(\iota(\prod_{p=1}^{\infty} \C[-2p]))$. Consequently
$$\mathbf{F}\simeq \Sigma^{-1} \ku\C/\Z\oplus \Sigma^{-1} H(\iota(\prod_{p=1}^{\infty} \C[-2p]))\ .$$
We have a long exact sequence
$$\dots \to \pi_{d+1}(\DD^{-}(X))\to \mathbf{F}^{-d}(X)\to \hat{ku}^{-d}(X)\stackrel{R}{\to}\pi_{d}(\DD^{-}(X))\to \dots \ .$$ 
Note that
$$H^{-k}(DD^{-}(X))\cong \prod_{p\in \Z} H^{-k}(\sigma^{\ge p}\Omega(X)[2p])\cong \prod_{p\in \Z} H^{2p-k}(\sigma^{\ge p}\Omega(X))\ .$$
If this group  is non-zero, then  for some $p\in \Z$ we have the relations
$$2p-k\le \dim(X)\ , \quad 2p-k\ge p\ .$$
This implies $\dim(X)\ge k$. Hence for $\dim(X)\le d$ we have $$\pi_{d+1}(\DD^{-}(X))=H^{-d-1}(DD^{-}(X))=0\  .$$  
Therefore
$\widehat{ku}^{-d}_{flat}(X)\cong \mathbf{F}^{-d}(X)$.
Since also
$$  \Sigma^{-1} H(\iota(\prod_{p=1}^{\infty} \C[-2p]))^{-d}(X  )=\prod_{p=1}^{\infty}H^{-2p-d-1}(X;\C)=0$$
we have the isomorphism
$$\mathbf{F}^{-d}(X)\cong \ku\C/\Z^{-d-1}(X)\ .$$
\hB

We decompose 
the curvature morphism $R$ in \eqref{sjaskjd8jkhdkqwdqw}   into a product of  components
$R(p)$, $p\in \Z$, according to the  product decomposition of $ DD^{-}$,  see    \eqref{dhwqdhqwkdqwdqwd}.
Similary we decompose $\reg_{X}$ into a product of components
$\reg_{X}(p)$ for $p\in \Z$. Note that $\pi_{d}(\DD^{-}(d)(X))\cong \Omega^{d}_{cl}(X)$.

\begin{lem}\label{jjkhdkjqhwdkjwqdwqdwqdwqd} 
 If    $\dim(X)\le  d$  and $x\in   K_{d} (C^{\infty}(X))$, then we have
$$R(p)(\hat\reg_{X}(x) )=\left\{\begin{array}{cc}0&p\not=d\\[0.3cm]\reg_{X}(d)(x)\in \Omega_{cl}^{d}(X)&p=d
\end{array}\right\}\ .$$
 \end{lem}
 \proof
Let $x\in   K_{d}(C^{\infty}(X))) $. 
 For $p\in \Z$ the component $R(p)(x)$ of the curvature of $x$ is represented by  the class  $\reg_{X}(p)( x ) \in H^{-d}(DD^{-}(p)(X))$.
By the calculations in the proof of Lemma \ref{kjewklfjwelkfewfwef} it can only be non-trivial if $p=d=\dim(X)$. Note that
$H^{-d}(DD^{-}(d)(X))\cong \Omega_{cl}^{d}(X)$.
 \hB
 
%
%
%
\begin{kor}\label{fjewfewfewfewfw}
If  $\dim(X)\le d-1$, then $\hat \reg_{X}$ maps to the subgroup of flat classes.
\end{kor}

In view of Lemma \ref{kjewklfjwelkfewfwef}  we can make the following definition.


%
\begin{ddd}\label{kdmlkqmdlqwdqwdwqd}
If  $\dim(X)\le d-1$, then we define $\sigma_{d}$ as the composition
\begin{equation}\label{dffdfdfdfdfrererer}\sigma_{d}:K_{d}(C^{\infty}(X)) \stackrel{\hat \reg_{X}}{\to} \widehat{ku}^{  -d}_{flat}(X)   \cong \ku\C/\Z^{-d-1}(X)\ .\end{equation}
 \end{ddd}

This finishes the proof of Theorem \ref{djekldjlwedwedewdew44dewde}. \hB

\begin{rem}{\rm 
At the moment we have no example which shows that the map $\sigma_{d}$ can be non-trivial if $\dim(X)<d-1$.
We refer to Example \ref{dhqwkdwqkdqwdqwdqwdqwd}, 7. and 8. for some vanishing results in this direction.}
\end{rem}

\subsection{Restriction to relative $K$-theory and explicit calculations}


In this subsection we derive an explicit formula for the restriction of $\sigma_{d}$ defined in \ref{kdmlkqmdlqwdqwdwqd} to topologically trivial classes. The result will be formulated as Corollary \ref{hqkdjhqwkjdhqwd}. It will be used in the proof of Theorem  \ref{dlkdjqwkdjwqlkdqwdqwd}.
The idea for the calculation of $\sigma_{d}(x)$ for a topologically trivial class $x\in K_{d}(X^{\infty}(X))$ is to use   a deformation $\tilde x$ of $x$ to zero. Using the homotopy formula  for  $\hat \ku$ the class $\sigma_{d}(x)$ can be expressed in terms of a transgression of the curvature of $\hat \reg_{X}(\tilde x)$.

\bigskip


 Let  $X$  be a manifold.
From the definition of $\hat \ku$ and the differential regulator $\hat \reg_{X}$ we get the diagram 

\begin{equation}\label{sjaskjd8jkhdkqwdqw1} \xymatrix{\bK^{rel}_{C^{\infty}(X)}\ar@{..>}[r]^{ \reg^{rel}_{X}}\ar[d]^{\partial}&i_{X}^{*}\Sigma\DD\ar@{=}[r]\ar[d]^{a }&i_{X}^{*}\Sigma\DD\ar[d]\\ \bK_{C^{\infty}(X)}\ar[r]^{\hat \reg_{X}}\ar[d]&i_{X}^{*}\hat{\mathbf{ku}} \ar[r]^{R}\ar[d]^{I}&  i^{*}_{X}\DD^{-} \ar[d]\\\bK^{top}_{C^{\infty}(X)}\ar[r]&i_{X}^{*}\underline{\mathbf{ku}}\ar[r]^{\ch^{cw}}& i_{X}^{*} \DD^{per}}\ .\end{equation} 
 The upper right square  follows from the fibre sequence \eqref{kdkjdkqwjdklwjqdlkqwdqwdwqd} and the fact that the lower right square is a pull-back, namely
  the definition \eqref{sjaskjd8jkhdkqwdqw} of $\hat \ku$.
  The columns are fibre sequences. 
  
  \begin{rem}{\rm  If $X$ is a compact manifold, then the columns in \eqref{sjaskjd8jkhdkqwdqw1}
   are instances of the homotopification  sequence \eqref{jjkfewfwefewfewffewf4333}. Indeed, the whole diagram can then be obtained by applying the homotopification sequence  to the middle raw. See Example \ref{kjdqlwkdlqwdqwldlqwd}.
 }\hB\end{rem}

\begin{kor}\label{hqkdjhqwkjdhqwd}
We have the equality
$$a\circ   \reg^{rel}_{X}=\hat \reg_{X}\circ \partial:K_{d}(C^{\infty}(X))\to \widehat{ku}^{-d}(X)\ .$$
\end{kor}

\bigskip

\bigskip

We assume that $X$ is a manifold of 
dimension $\dim(X)=d-1$. 
We use the homotopy formula for $\hat{\ku}$ in order to provide  a   
formula for $\sigma_{d}(x)$ for certain topologically trivial classes $x\in K_{d}(C^{\infty}(X))$.
 We consider a class $\tilde x\in K_{d}(C^{\infty}(I\times X))$ which has the property that
$\tilde x_{|\{0\}\times X}=0$. We define $$x:=\tilde x_{|\{1\}\times X}\in K_{d}(C^{\infty}(X))\ .$$
By construction the class  $x$ is topologically trivial, i.e. it belongs to the kernel of $K_{d}(C^{\infty}(X))\to K_{d}^{top}(C^{\infty}(X))$.

Since the degree of $\tilde x$ and the dimension of $I\times X$ match,
by Lemma \ref{jjkhdkjqhwdkjwqdwqdwqdwqd} the only non-trivial component of the regulator of $\tilde x$ is $$\reg_{I\times X}(d)(\tilde x)\in \Omega^{d}_{cl}(I\times X)\ .$$ Furthermore we have a map
$$i_{d}:\Omega^{d-1}_{cl}(X)\to H^{d-1}(\sigma^{<d}\Omega(X))\to H^{-d-1}((\sigma^{<d}\Omega(X))[2d])\to  \pi_{d}( \Sigma \DD (X))\ .$$  
\begin{prop}\label{dwjdwqdjwdlwkdqwd}  We have 
$$\hat \reg_{X}(x)= a (i_{d}\int_{I}  \reg_{I\times X}(d)(\tilde x))\ .$$
\end{prop}
\proof
By the homotopy formula \cite[(27)]{2013arXiv1311.3188B} we have 
$$\hat \reg_{X}(x)=a (\int_{I} R  (\hat \reg_{I\times X}(\tilde x)))\ .$$
By Lemma \ref{jjkhdkjqhwdkjwqdwqdwqdwqd} we have $$R (d)(\hat \reg_{I\times X}(\hat x))= \reg_{I\times X}(d)(\tilde x)\ ,$$ and all other factors of the curvature vanish.
So our final formula is
$$\hat \reg_{X} (x)=a ( i_{d}\int_{I}  \reg_{I\times X}(d)(\tilde x))\ .$$
\hB
 
\begin{ex}\label{dhqwkdwqkdqwdqwdqwdqwd}{\rm 
We now perform some explicit calculations:
\begin{enumerate}
\item We consider $X=S^{1}$. Let $u:S^{1}\to \C$ be the inclusion. We consider $u$ as a unit in $C^{\infty}(S^{1})$ and $\iota(u)\in K_{1}(C^{\infty}(S^{1}))$. Note that
$$\pi_{1}(\DD^{-}(S^{1}))\cong \pi_{1}(\DD^{-}(1)(S^{1}))\cong \Omega^{1}(S^{1})\ .$$ Under this identification 
 have $$\reg_{S^{1}}(1)(\iota(u))= \frac{1}{2\pi i}\frac{du}{u}\ .$$
\item
Next we consider the inclusion $t:\R_{+}\to \C$ as a unit in $C^{\infty}(\R_{+})$. 
We obtain $$\reg_{\R_{+}}(1)(\iota(t))= \frac{1}{2\pi i}\frac{dt}{t}\ .$$
\item
We now consider the inclusion $z:\C^{*}\to \C$. We write
$z=u(z) t(z)$ with $u:=z/|z|^{-1}$ and $t:=|z|$.
Then we get by naturality and additivity of the regulator
$$\reg_{\C^{*}}(1)(\iota(z))= \frac{1}{2\pi i}\left(\frac{d(z/|z|^{-1})}{z/|z|^{-1}}+\frac{d|z|}{|z|}\right)= \frac{1}{2\pi i}\frac{dz}{z}\ .$$
\item
We  consider the unit  $\exp(z)=\exp^{*}z\in C^{\infty}(\C)$. We get
$$\reg_{\C}(1)(\iota(\exp(z))=\exp^{*} (\reg_{\C^{*}}(1)(\iota(z)))= \frac{1}{2\pi i} dz\ .$$
\item
Let now $X$ be a smooth manifold and $f\in C^{\infty}(X)$. Then
we have \begin{equation}\label{hgdhjgqwjdqwdqwdwqdqwdwqd88}\reg_{X}(1)(\iota(\exp(f))=f^{*}(\reg_{\C}(1)(\iota(\exp(z)))= \frac{1}{2\pi i}df\ .\end{equation}
\item
The regulator $\reg_{X}:K_{*}(C^{\infty}(X))\to \pi_{*}(\DD^{-}(X))$ is
given by the composition of the Goodwillie-Jones Chern character and the map \eqref{kkckaslckjaslkcjaslkcjasc}. Since $C^{\infty}(X)$ is commutative  both maps are in fact multiplicative, where the ring structure on $\DD^{-}$ is induced by the obvious bigraded differential algebra structure on $\prod_{p\in \Z}(\sigma^{\ge p}\Omega)[2p]$.

Now assume that    $\dim(X)=d-1$. We     consider a collection $u_{2},\dots,u_{d}$ of units in $C^{\infty}(X)$.
Let $t\in I$ be the coordinate. On $I\times X$ we consider the collection of $d$ units
$$\exp(tf),u_{2},\dots,u_{d}$$
and define
$$\tilde x:=\iota(\exp(tf))\cup \iota(u_{2})\cup \dots\cup \iota(u_{d})\in K_{d}(C^{\infty}(I\times X))\ .$$

By multiplicativity of the regulator and \eqref{hgdhjgqwjdqwdqwdwqdqwdwqd88} we have
$$\reg_{I\times X}(d)(\iota(\exp(tf))\cup \iota(u_{2})\cup \dots\cup \iota(u_{d})) = \frac{1}{(2\pi i)^{d}} (f dt+ tdf) \wedge \frac{du_{2}}{u_{2}}\wedge \dots \wedge\frac{du_{d}}{u_{d}}$$
(and this is the only non-trivial component).
From Proposition \ref{dwjdwqdjwdlwkdqwd} we conclude that
\begin{equation}\label{hkhkjwdhkwjqdhjkqwdwqdqwd}\hat \reg_{X}(\iota(\exp(f))\cup \iota(u_{2})\cup \dots\cup \iota(u_{d}))=a ( i_{d}( \frac{1}{(2\pi i)^{d}}f\frac{du_{2}}{u_{2}}\wedge \dots \wedge\frac{du_{d}}{u_{d}}))\ .\end{equation}
\item If we consider  in 6.  more than $d$ factors, then we get a trivial result.
Assume that $\dim(X)=d-1$ and $ u_{2},\dots,u_{d+1}$ are invertible  smooth functions, and $f$ is a smooth function. Then we have
 $$\hat \reg_{X}\left(\iota(\exp(f))\cup \iota(u_{2})\cup \dots \cup\iota(u_{d+1})\right)=0\ .$$
To this end, using multiplicativity of $\hat \reg_{X}$, we write the left-hand side as a product in $\hat \ku$-theory $$\hat \reg_{X}\left(\iota(\exp(f))\cup \iota(u_{2})\cup \dots \cup\iota(u_{d})\right) \cup \hat \reg_{X}(\iota(u_{d+1}))\ .$$ 
Using \eqref{hkhkjwdhkwjqdhjkqwdwqdqwd} and the rule $a(\omega)\cup x=a(\omega\wedge R(x))$ for the product in $\hat \ku$ we rewrite this as
$$a ( i_{d}( \frac{1}{(2\pi i)^{d}}f\frac{du_{2}}{u_{2}}\wedge \dots \wedge\frac{du_{d}}{u_{d}}))\cup  \hat \reg_{X}(\iota(u_{d+1}))=  a ( i_{d+1}( \frac{1}{(2\pi i)^{d+1}}f\frac{du_{2}}{u_{2}}\wedge \dots \wedge\frac{du_{d+1}}{u_{d+1}}))
 \ .$$
 The argument of $a$ on the right-hand side is a $d$-form on a manifold of dimension $d-1$ and therefore vanishes.
\item
Let again $\dim(X)\le d-1$, $x\in K_{d}(C^{\infty}(X))$ and $y\in K_{d^{\prime}}(C^{\infty}(X))$ for some $d^{\prime}\in \nat$.
Then
$$\hat \reg_{X}(x\cup y)=\hat \reg_{X}(x)\cup \hat \reg_{X}(y)=\sigma_{d}(x)\cup  y^{top}\ ,$$
where $y^{top}:=I(\hat \reg_{X}(y))\in \ku^{-d^{\prime}}(X)$, using the $\ku$-module structure of $\ku\C/\Z$ on the right-hand side.

Note that $\ku(X)$ has a multiplicative Atiyah-Hirzebruch filtration $(F^{p}\ku(X))_{p\in \Z}$, where  $F^{p}\ku(X)$ consists   classes whose pull-backs to $p-1$-dimensional $CW$-complexes vanish. 
For an invertible function $u$ we have $\iota(u)^{top}\in F^{1}\ku^{-1}(X)$. 
If $y=\iota(u_{1})\cup \dots \cup \iota(u_{d})$, then $y^{top}=0$. We conclude that
$$\hat \reg_{X}\left( x\cup  \iota(u_{1})\cup \dots \cup \iota(u_{d}) \right)=0\ .$$

\end{enumerate}
}\hB\end{ex}

\section{An index theorem}\label{kwndqlwdqwdqwdwqdwqdwdwqd}

\subsection{The index pairing}

In this subsection we introduce the pairing between Dirac operators and the Hopkins-Singer version of  differential periodic complex $K$-theory. This pairing refines the pairing between periodic complex $K$-theory and the $K$-homology. 
If the Dirac operator comes from a $Spin^{c}$-structure, then our  pairing is a special case of the integration in differential cohomology. In any case,  the pairing  has a simple description in terms of 
standard constructions of local index theory.

\bigskip

First we introduce the Hopkins-Singer version $\widehat{KU}^{*}(-)$ of differential periodic complex $K$-theory $\bku$.  By   $$\ch^{cw}_{per}:\bku\to H(\iota(\prod_{p\in \Z}\C[2p]))$$ we denote  the usual Chern character. We will use the same symbol also in order to denote  the composition of morphisms of sheaves of spectra
$$\ch^{cw}_{per}:\underline{\bku}\stackrel{\ch^{cw}_{per}}{\to} \underline{H(\iota(\prod_{p\in \Z}\C[2p]))}\stackrel{\eqref{jskljkjlkjwsqws67567572613}}{\simeq} \DD^{per}\ .$$
\begin{rem}{\rm  The Chern character  $\ch^{cw}$ discussed in Lemma \ref{jjljldkqjwdqwdqwd} is equivalent to the composition  
$$\ku\to \bku\stackrel{\ch_{per}^{cw}}{\to} H(\iota(\prod_{p\in \Z}\C[2p]))\ ,$$
where the first map 
\begin{equation}\label{egdjhdgjhd}\ku\to \bku\end{equation}is the connective covering map of $\bku$.

}\hB\end{rem}
We set
$${}^{\ge 0}\DD^{per}:=H(\iota(\sigma^{\ge 0} DD^{per}))\ .$$ The inclusion of sheaves of chain complexes 
$\sigma^{\ge 0} DD^{per}\to   DD^{per}$ induces a morphism of sheaves of spectra
${}^{\ge 0}\DD^{per}\to \DD^{per}$.

\begin{ddd}\label{jdkjwldqwdqwdqwdqwd}  We define the Hopkins-Singer differential cohomology theory associated to $\bku$  as the sheaf of spectra $\hat \bku\in \Sh_{\Sp}(\Mf)$ given by the 
  following pull-back   $$\xymatrix{\hat \bku\ar[r]^{R}\ar[d]^{I}&{}^{\ge 0}\DD^{per}\ar[d]\\
\underline{\bku}\ar[r]^{\ch^{cw}_{per}}&\DD^{per}}\ .$$
We further define the differential periodic complex $K$-theory  groups of a manifold $X$  by  $$\widehat{KU}^{*}(X):=\pi_{-*}(\hat \bku(X))$$ (compare with Definition \ref{fljfkelwjflewfewfef}).
\end{ddd}

The differential  periodic complex $K$-theory in degree zero fits into the exact sequence \begin{equation}\label{sjxlkasjxklasjcxsaxjsxjsax}\bku^{-1}(X)\stackrel{\ch^{cw}_{per}}{\to}  DD^{per}(X)^{-1}/\im(d)\stackrel{a_{\hat \bku}}{\to} \widehat{KU}^{0}(X)\stackrel{I}{\to} \bku^{0}(X)\to 0\ .\end{equation}
 This exact sequence is one of the basic features of differential  periodic complex $K$-theory, see e.g.   \cite[Prop. 2.20]{MR2664467}.

By \cite{MR2608479} for a compact manifold $X$  the group $\widehat{KU}^{0}(X)$
is canonically isomorphic to the differential $K$-theory groups defined using  geometric models  \cite{MR2664467}, \cite{MR2732065}. In these models  geometric vector bundles are cycles for differential $K$-theory  classes. Recall that
a 
geometric $\Z/2\Z$-graded vector bundle $\bV:=(V,h,\nabla)$ is a triple consisting of a $\Z/2\Z$-graded complex vector bundle $V\to X$, a hermitean metric $h$ on $V$ such that the even and odd summands are orthogonal, and a connection $\nabla$ which preserves $h$ and the grading.
In the geometric models for  $\widehat{KU}^{0}(X)$ a
geometric $\Z/2\Z$-graded vector bundle $\bV:=(V,h,\nabla)$  tautologically represents   a class
$$[\bV]\in \widehat{KU}^{0}(X)\ .$$
We refer to \cite[Sec. 6.1]{2013arXiv1311.3188B} for an alternative construction of this class 
in terms of a  cycle map.

\bigskip

 We now assume that $X$ is a  closed Riemannian manifold of odd dimension $d$. 
 We further assume that we are given a generalized Dirac operator  $\Dirac$ on $X$.
 By definition, $\Dirac$ is the Dirac operator associated to a Dirac bundle, see e.g. \cite[Sec. 3.1]{MR2191484}.

 \begin{rem}{\rm A generalized Dirac operator provides a $K$-homology class  which
 can be paired with  $K$-theory classes on $X$. The basic idea of the following Lemma is that the Dirac operator as a geometric object provides a sort of differential refinement of its $K$-homology class which can be paired with differential $K$-theory classes.  The map $\rho_{\Dirac}$ defined  in Proposition \ref{hdkjhekwdewdewdewd} below only captures the secondary information contained in this pairing. Its value on a differential $K$-theory class can be considered as the reduced $\eta$-invariant of the Dirac operator twisted with this class.
A very similar construction has been used in order to define the intrinsic universal $\eta$ invariant in  
 \cite{2011arXiv1103.4217B}.}\hB
 \end{rem}
\begin{prop} \label{hdkjhekwdewdewdewd}
 We have a canonical evaluation map
 $$\rho_{\Dirac}:\widehat{KU}^{0}(X)\to \C/\Z\ .$$
 \end{prop}
 \proof
  Let $x\in \widehat{KU}^{0}(X)$.  In view of the sequence \eqref{sjxlkasjxklasjcxsaxjsxjsax} we can choose   a geometric $\Z/2\Z$-graded vector bundle $\bV:=(V,h,\nabla)$ and a form 
  $\gamma\in DD^{per}(X)^{-1}/\im(d)$ such that  the following identity holds true in $\widehat{KU}^{0}(X)$: $$x=[\bV]+a_{\hat \bku}(\gamma)\ .$$ 
 We need the following standard  constructions from local index theory:
 \begin{enumerate}
 \item We form the twist $\Dirac\otimes \bV$  of the Dirac operator by $\bV$ (see \cite[Sec. 3.1]{MR2191484} for details if necessary). 
 \item 
 We let $\xi(\Dirac\otimes \bV)\in \R/\Z$ denote the reduced $\eta$-invariant  of $\Dirac\otimes \bV$ given by
 \begin{equation}\label{gdghgjwedewd7787864324324}\xi(\Dirac\otimes \bV):=[\frac{\eta(\Dirac\otimes \bV)+\dim(\ker(\Dirac\otimes \bV))}{2}]\ ,\end{equation}
 where $\eta(\Dirac\otimes \bV)$ is the  Atiyah-Patodi-Singer $\eta$-invariant introduced  in \cite{MR0397797}.
 \item 
 We let $\hA(\Dirac)\in \prod_{p\in \Z} (\Omega (X,\Lambda)[2p])^{0}_{cl}  $ denote the local index form associated to $\Dirac$ where $\Lambda$ denotes the orientation twist of $X$.  
 \begin{rem}\label{ewjflewjflewfewfewf}{\rm The local index form  has the following  explicit description.  Locally on $X$ we can write $\Dirac=\Dirac_{spin}\otimes \bE$ for the spin Dirac operator $\Dirac_{spin}$ and a geometric $\Z/2\Z$-graded twisting bundle $\bE=(E,h^{E},\nabla^{E})$.  If we can write $\Dirac$ in this way, then
  $$\hA(\Dirac)=([\hA(\nabla^{LC})\wedge \ch(\nabla^{E})]_{2p})_{p\in \Z}\in \prod_{p\in \Z}(\Omega (X,\Lambda)[2p])^{0}_{cl}\ ,$$
 where $\nabla^{LC}$ is the Levi-Civita connection of $X$, 
$[\omega]_{2p}$ denotes the degree-$2p$-component of the inhomogeneous even form $\omega$, and
$\hA(\nabla)$ and $\ch(\nabla)$ are the usual characteristic forms defined in terms of the curvature of
the connections (including the $2\pi i$-factors), see \cite[Sec. 4.3]{MR2191484} for explicit formulas.

The following observation will make it unnecessary to use the explicit formula for the index density. 
This fact will be particularly helpful in the proof of Lemma \ref{djqwjnkddqwdqwd} below.
We define
the integral
\begin{equation}\label{dkqwdkwqdlkdjwqkldjwqlkdj}\int_{X}: \prod_{p\in \Z} \Omega(X,\Lambda)[2p] \to \C\ ,  \quad \int_{X}(\omega(p))_{p\in \Z}:=\sum_{p\in \Z} \int_{X }[\omega(p)]_{\dim(X)}\ .\end{equation}
It induces an evaluation of cohomology classes which we will denote by the same symbol. 
By the Atiyah-Singer  index theorem we can calculate for every class  $u\in \bku^{-1}(X)$ the  index pairing
by
\begin{equation}\label{kdklqjwdljwqdlwqdwqdqwdqwdwd}\langle [u],[\Dirac]\rangle=\int_{X} [\hA(\Dirac)]\cup  \iota_{d+1} \ch^{cw}_{per}(u)\in \Z\ ,\end{equation}
where $\iota_{d+1}$ is the shift isomorphism defined in \eqref{shiftiso}.


 }\hB\end{rem}
 \end{enumerate}



 Using the integral \eqref{dkqwdkwqdlkdjwqkldjwqlkdj} we   now define
\begin{equation}\label{edewdwede}\rho_{\Dirac}(x):=\xi(\Dirac\otimes \bV)+[\int_{X} \hA(\Dirac)\wedge \iota_{d+1}\gamma]\in \C/\Z\ .\end{equation}  

We must show that
  $\rho_{\Dirac}$ is a well-defined homomorphism.

  \begin{enumerate}
  \item First observe that  the right-hand side  of \eqref{edewdwede}   does not depend on the choice of $\gamma$.
 Indeed, by \eqref{sjxlkasjxklasjcxsaxjsxjsax} two choices differ by a closed form representing an element in the image of $\ch^{cw}_{per}:\bku^{-1}(X)\to DD^{per}(X)^{-1}/\im(d)$, and the integral of the product of those elements with $\hA(\Dirac)$ is an integer by the Atiyah-Singer index theorem, see Remark \ref{ewjflewjflewfewfewf}.
 \item 
 We observe that the right-hand side of \eqref{edewdwede} is invariant under stabilization by bundles which admit an odd $\Z/2\Z$ symmetry. 
\item  Next we observe, using the variation formula
for the classes $[\bV]$ (the homotopy formula for $\widehat{KU}^{0}$) and $\xi(\Dirac\otimes \bV)$, that
the right-hand side of \eqref{edewdwede} does not depend on the choice of the geometry of $V$.
\item
If we choose a different bundle $\bV^{\prime}$ and form $\gamma^{\prime}$,
then after stabilization we can assume that
$V\cong V^{\prime}$ as graded bundles. Therefore the right-hand side of \eqref{edewdwede} does not depend on the choice of $\bV$.
\item Finally we observe that $\rho_{\Dirac}$ is a homomorphism. 
\end{enumerate}
\hB

%

\begin{rem}\label{djlwjdqwdqwdqwdqwdqwdwqd}{\rm
For any integer $k\in \Z$ we can define a version $\widehat{KU}^{k,*}$ of differential $K$-theory by replacing the cut-off ${}^{\ge 0}$   in Definition \ref{jdkjwldqwdqwdqwdqwd} with ${}^{\ge k}$. 
Assume that $X$ is a Riemannian spin manifold of odd dimension $d$. Then the map $p:X\to *$ is differentially $K$-oriented and we have an integration
$$\hat p_{!}:\widehat{KU}^{0,0}(X)\to \widehat{KU}^{-d,-d}(*)\cong \C/\Z\ .$$
We refer to  \cite[Sec. 3]{MR2664467} and \cite[Sec. 4.10 and 4.11]{skript} for details on the integration in differential cohomology.

Let us now assume that $\Dirac=\Dirac_{spin}\otimes \bE$  for some twist $\bE=(E,h^{E},\nabla^{E})$.
From \cite[Cor. 5.5]{MR2664467} we conclude that $\rho_{\Dirac}$ can be expressed in terms of the integration $ \hat p_{!}$ in differential $K$-theory as follows
$$\rho_{\Dirac}(x)=\hat p_{!}([\bE]\cup x)\ , \quad x\in \widehat{KU}^{0}(X)\ .$$ 

The spin structure on $X$ provides the underlying topological $K$-orientation of $X$ given by the fundamental class $[\Dirac_{spin}]\in \bku_{d}(X)$. The restriction of $\rho_{\Dirac}$ to the flat subgroup corresponds under the identification 
\begin{equation}\label{jkhfjkehfehewfwefe}\widehat{KU}^{0}_{flat}(X)\cong \bku\C/\Z^{-1}(X)\end{equation}
(this is the analog of Lemma \ref{kjewklfjwelkfewfwef})
to the evaluation pairing
$$\langle -\cup [E] ,[\Dirac_{spin}]\rangle:\bku\C/\Z^{-1}(X)\stackrel{-\cup [E]}{\to} \bku\C/\Z^{-1}(X) \stackrel{\langle-,[\Dirac_{spin}]\rangle}{\to} \bku\C/\Z^{-d-1}(*)\cong \C/\Z\ ,$$
where the first map uses the $\bku$-module structure of $\bku\C/\Z$.

%
%
%
 }\hB \end{rem}

\begin{rem}{\rm
In this remark we explain the relation between the evaluation map $\rho_{\Dirac}$ and the index theorem for flat vector bundles by Atiyah-Patodi-Singer \cite[Thm. 5.3]{MR0397799}.
Let $\bV$ be a $\Z/2\Z$-graded flat geometric bundle of virtual dimension zero. It represents a class
$[\bV]\in \widehat{KU}^{0}_{flat}(X)\cong \bku\C/\Z^{-1}(X)$.  
In this case
$$\rho_{\Dirac}([\bV])=\xi(\Dirac\otimes \bV)$$ is exactly the analytic index of the flat bundle introduced by Atiyah-Patodi-Singer.
Their index theorem for flat bundles states that this analytic index is equal to
the pairing of the class $[\bV]$ with the $K$-homology class of $\Dirac$. This also follows from the last assertion in Remark \ref{djlwjdqwdqwdqwdqwdqwdwqd}.
 } \hB
\end{rem}


 \begin{ex}\label{jckbcwkecewce}{\rm
 We consider $S^{1}$ as a spin manifold with the standard metric of length $1$ and  with the non-bounding spin structure. The spinor bundle is one-dimensional and can be trivialized such that $\Dirac_{spin}=i\partial_{t}$. Assume now that $\bL$ is a geometric line bundle with  holonomy $v\in U(1)$.
Then we can trivialize $\bL$ such that its connection is given by $\nabla^{L}=d-\log(v) dt$.
We get  $$\Dirac_{spin}\otimes \bL=i(\partial_{t}-\log(v))\ .$$  Its spectrum is
$\{2\pi  n -\log v\:|\: n\in \Z\}$
with multiplicity $1$. For $v\not=1$
we get by an explicit calculation 
$$\eta(\Dirac_{spin}\otimes \bL)=1-\frac{\log(v)}{\pi i}\ ,$$ where the branch of the logarithm is chosen such that $\frac{\log(v)}{\pi i}\in (0,2)$. 
Using \eqref{gdghgjwedewd7787864324324}  and \eqref{edewdwede} we get
$$\rho_{\Dirac_{spin}}([\bL])=[\frac{1}{2}-\frac{\log(v)}{2\pi i}]_{\C/\Z}$$
because in this case we can take $\gamma=0$. 
This formula holds true also for $v=1$.

\bigskip

If $\bL$ is trivial, then $[\bL]=1$ and   we have  $\rho_{\Dirac_{spin}}(1)=[\frac{1}{2}]$.
We have an isomorphism $\hat \bku^{0,0}_{flat}(S^{1})\cong \C/\Z$ which maps
$[\bL]-1$ to $\frac{\log v}{2\pi i}$. With this identification  the  restriction of the evaluation map to the flat subgroup  is given by the homomorphism
$$\rho_{\Dirac_{spin}}:\C/\Z\to \C/\Z\ , \quad [z]\mapsto [-z]\ .$$
  }
 \hB\end{ex}


 \subsection{Fredholm modules}\label{dhwiqdjoqwidqwdqwdqwdwqdqdqwd123}
 
 We consider a closed  Riemannian manifold $X$ of odd dimension $d$ with a generalized Dirac operator $\Dirac$ associated to a Dirac bundle $\bE$. 
The Dirac operator $\Dirac$ gives rise to a $d+1$-summable Fredholm module $(H,P )$ over $C^{\infty}(X)$ as follows (see \cite{MR823176}):  
\begin{enumerate}
\item The Hilbert space of the Fredholm module is
$H:=L^{2}(X,E)$. The   algebra $C^{\infty}(X)$ acts on $H$ in the usual way by multiplication operators. \item The operator $P\in B(H)$ is the orthogonal projection $P^{+}$ onto the positive   eigenspace of $ \Dirac$
 \item The condition that $(H,P)$ is   $d+1$-summable means that $$ [P,f]\in \cL^{d+1}(H)\ ,$$
for all $f\in C^{\infty}(X)$, where $\cL^{d+1}(H)$ denotes the $d+1$'th Schatten class.
\end{enumerate}

\begin{rem}{\rm
In some references odd Fredholm modules are denoted by $(H,F)$, where $F\in B(H)$
is a selfadjoint involution such that $[F,f]\in \cL^{d+1}(H)$. 
The relation with our notation is given by the equation $F=P-(1-P)$.
}
\end{rem}


We let $\cM^{d}$ be the universal algebra for
 $d+1$-summable Fredholm modules  introduced by Connes-Karoubi \cite{MR972606}. Then we get a homomorphism
 \begin{equation}\label{fevgerhjvhjbjhervvev}b_{\Dirac}:C^{\infty}(X)\to \cM^{d}\end{equation} classifying the Fredholm module  $(H,P)$. Note that $b_{\Dirac}$ is uniquely determined up to unitary equivalence.

 \begin{rem}\label{kljkljddlqwdqwdqwd}
{\rm  In this remark we give an explicit description of $b_{\Dirac}$. The algebra $\cM^{d}$ for odd $d\in \nat$ is  a subalgebra of the algebra of  $2\times 2$-matrices  of bounded   operators on the standard separable Hilbert space $\ell^{2}$
 consisting of the matrices    \begin{equation}\label{}\left(\begin{array}{cc}a_{11}&a_{12}\\a_{21}&a_{22}\end{array}\right) \ , \quad a_{12}, a_{21}\in \cL^{d+1}(\ell^{2})\ , \quad  a_{11}, a_{22}\in B(\ell^{2})\ .\end{equation}   

Let $P^{+},P^{-}$ be the positive and non-positive  spectral projections of $\Dirac$.
Then we choose identifications $\ell^{2}\cong \im(P^{+})\cong \im(P^{-})$.
The homomorphism $b_{\Dirac}:C^{\infty}(X)\to \cM^{d}$ is then given by
\begin{equation}\label{dqwdwqdwqdqwdqwdw324324234}f\mapsto \left(\begin{array}{cc}P^{+}fP^{+}&P^{+}fP^{-}\\P^{-}fP^{+}&P^{-}fP^{-}\end{array}\right)\ .\end{equation}
}\hB\end{rem}

\subsection{The multiplicative character} \label{lkewdjlkwefwfwfewffefewfewf}
 
In \cite[4.10]{MR972606} Connes and Karoubi  constructed the "multiplicative" character
\begin{equation}\label{edwededwedwe32423423554345}\delta:K_{d+1}(\cM^{d})\to \C/\Z\ .\end{equation}
 In this subsection we explain how the construction of the multiplicative character  $\delta$ fits into the framework of differential cohomology. The details of the construction will be needed later in Subsection \ref{djoiqwjdoqwdjwqdwqdd}.
 
 \bigskip
 
 We consider a unital locally convex algebra $A$. It has a natural diffeological structure $A^{\infty}$, see Example \ref{dhqwkdhwqkdhqkjdqwdwq}, 5.  From the sheaf of algebras $A^{\infty}$ we derive  the    sheaves  of spectra $$\CC^{-}_{A}:=L(H(\iota(CC^{-}(A^{\infty}))))\ , \quad \bK_{A}\stackrel{\eqref{deldqlkdwdwqddjh9879879}}{=}L(\bK(A^{\infty}))\ .$$


  The  Goodwillie-Jones Chern character \eqref{djkejdkejdewd777wed} gives a morphism 
\begin{equation}\label{smkjlkwjsljwqswqs}\ch^{gj}:\bK_{A}\to \CC^{-}_{A}\ .
\end{equation}

\begin{rem}{\rm
In principle we  want to apply the homotopification sequence \eqref{jjkfewfwefewfewffewf4333} to $\ch^{gj}$. This leads to the problem of understanding the homotopification of $\CC^{-}_{A}$. The known facts are contained e.g. in \cite[Cor. 4.1.2]{MR2409415}. In particular the homotopification of $\CC^{-}_{A}$ is equivalent to the homotopification of its periodic version $\CC^{per}_{A}$, and the homotopification of the cyclic homology  
$\CC_{A}$ vanishes. The problem is that $\CC^{per}_{A}$ is not known to be homotopy invariant.
We will get a better theory if we use the continuous versions of cyclic homology. The main advantage is that the continuous periodic cyclic homology for complete locally convex algebras is known to be diffeotopy invariant, see Theorem \ref{jhdkjwqdqwdqwdwdqwdqwdqwdwqdwqd}.
}\hB 
\end{rem}

If we define the cyclic bicomplex $\cB C^{cont}(A)$ of a locally convex algebra $A$ similarly as in \cite[5.1.7]{MR1600246} but using projective tensor products, then we get the
continuous versions of cyclic, negative cyclic and and periodic cyclic homology complexes \begin{equation}\label{soqwjsoqwjswqswqs}CC^{cont}(A)\ , \quad CC^{cont,-}(A)\ ,\quad CC^{cont,per}(A)\ .\end{equation}
In the natural extension of the  notation of \cite[5.1.7]{MR1600246} to the continuous case these complexes would have been denoted by
$\mathrm{Tot}\: \cB C^{cont}$, $\mathrm{ToT} \: \cB C^{cont,-}$, and 
$\mathrm{ToT}\: \cB C^{cont,per}$. 
We have an exact sequence of chain complexes
\begin{equation}\label{llwdlqwddqwdkllkklklklkll}0\to CC^{cont,-}(A) \to CC^{cont,per}(A)\stackrel{q}{\to} CC^{cont}(A)[2]\to 0\ .\end{equation}
Note that $A^{\infty}$ is a presheaf of locally convex algebras by Remark \ref{ldjlkqdjqwldqwdwqd7987}. 
In \eqref{llwdlqwddqwdkllkklklklkll} we can thus replace $A$ by $A^{\infty}$ in order to get an exact sequence of presheaves with values in $\Ch$.  
 Then we apply $L\circ H\circ \iota$ in order to get the fibre sequence of sheaves of spectra
\begin{equation}\label{dkjdwlkqjdlqwjdwqdwqdqwdwqd234}\Sigma \CC^{cont}_{A} \to \CC^{cont,-}_{A}\to \CC^{cont,per}_{A}\to \Sigma^{2}\CC^{cont}_{A} \end{equation}
which is very similar to \eqref{kdkjdkqwjdklwjqdlkqwdqwdwqd}.

We now use the well-known fact that the continuous periodic cyclic homology is diffeotopy invariant \cite[Theoreme 2.7]{MR1478702} (for Fr\'echet algebras) and \cite{cva} (for complete locally convex algebras):
\begin{theorem}\label{jhdkjwqdqwdqwdwdqwdqwdqwdwqdwqd}  
Assume that $A$ is a  complete locally convex algebra. Then
the projection $I\times M\to M$ induces a quasi-isomorphism
$$CC^{cont,per}( C^{\infty}( M,A))\to CC^{cont,per}(C^{\infty}(I\times M,A))\ .$$
\end{theorem}
 As a consequence, the sheaf $\CC^{cont,per}_{A}$ is homotopy invariant in the sense of  Definition \ref{kwhdwqjkdhkwdjwqddwqdwqdwqd32423423}. 

\bigskip 
From now on we assume that $A$ is complete.
We   apply the homotopification sequence \eqref{jjkfewfwefewfewffewf4333} to the morphism \eqref{smkjlkwjsljwqswqs}. Using the Definition \ref{dhqkdhwqkjdqwdqwdqwdwqd23424}
we obtain the upper two columns of the following diagram:
 
 \begin{equation}\label{slkkjle}\xymatrix{
\bK^{rel}_{A}\ar[r]\ar[d]^{\cA(\ch^{gj})}&\bK_{A}\ar[r]\ar[d]^{\ch^{gj}}&\bK^{top}_{A}\ar[d]^{\cH(\ch^{gj})}\ar[r]&\Sigma \bK^{rel}_{A}\ar[d]\\
\cA(\CC^{-}_{A})\ar[r]\ar@{-->}[d]^{ii}&\CC^{-}_{A}\ar[d]^{t}\ar[r]&\cH(\CC^{-}_{A})\ar[r]\ar@{..>}[d]^{i}&\cA(\CC^{-}_{A})\ar@{-->}[d]^{ii}\\
\Sigma \CC^{cont}_{A}\ar[r]&\CC^{cont,-}_{A}\ar[r]^{p}&\CC^{cont,per}_{A}\ar[r]&\Sigma^{2}\CC^{cont}_{A}}
\end{equation}
The lower sequence is \eqref{dkjdwlkqjdlqwjdwqdwqdqwdwqd234}. The map $t:\CC^{-}_{A}\to \CC^{cont,-}_{A}$ is induced by the canonical map from algebraic to projectively completed tensor products.
The composition $p\circ t$ maps to a homotopy invariant target. The dotted arrow marked by $i$ and the filler of the lower middle square are obtained from the universal property of the homotopification as a left adjoint in \eqref{dhqjdhkjqwdqwdw}.
The dashed arrows marked by $ii$ and the corresponding fillers are now induced naturally.

We now drop out the middle row and evaluate the diagram at $*$. We then get  the map of fibre sequences of spectra
 \begin{equation}\label{jkljdlqkwdqwdqwd}\xymatrix{
\bK^{rel}(A)\ar[r]\ar[d]^{\ch^{cont}_{rel}}&\bK(A)\ar[r]\ar[d]^{\ch^{cont}}&\bK^{top}(A)\ar[d]^{\ch^{cont}_{per}}\ar[r]&\Sigma \bK^{rel}(A)\ar[d]\\ 
\Sigma \CC^{cont}(A)\ar[r]&\CC^{cont,-}(A)\ar[r]^{p}&\CC^{cont,per}(A)\ar[r]^{q}&\Sigma^{2}\CC^{cont}(A)}
\end{equation}
which defines various versions of the continuous Chern character.

\begin{rem}{\rm
By \cite[Lemma 4.2.2]{MR2409415} the lower sequence in \eqref{jkljdlqkwdqwdqwd} can be identified with the homotopification
sequence of $\CC^{cont,-}(A)$. The whole diagram is thus the result of applying the homotopification sequence to the map  $t\circ \ch^{gj}:\bK(A)\to \CC^{cont,-}(A)$. The construction of the various versions of the continuous Chern characters  above is thus completely parallel to what is done in  \cite[Sec. 4.2]{MR2409415}. The diagram \eqref{jkljdlqkwdqwdqwd} is exactly the last diagram in  \cite[Sec. 4.2]{MR2409415}.
}\hB
\end{rem}

We finally  define the Chern character $\ch^{cont}_{top}$ by the following diagram which involves the factorization of $q$ over the $S$-operator:
\begin{equation}\label{jdklqwjdlkjqwdwd}\xymatrix{\bK^{top}(A)\ar[d]^{\ch^{cont}_{per}}\ar[rd]^{\ch^{cont}_{top}}\ar[rr]&&\Sigma \bK^{rel}(A)\ar[d]^{\ch^{cont}_{rel}}\\
\CC^{cont,per}(A)\ar@{..>}@/_1cm/[rr]^{q}\ar[r]&\CC^{cont}(A)\ar[r]^{S}&\Sigma^{2}\CC^{cont}(A)}
\end{equation}

For $\sharp\in \{\emptyset,-,per\}$ we let $HC^{cont,\sharp}_{*}(A):=H_{*}(CC^{cont,\sharp}(A))$ denote the respective versions  of continuous cyclic homology groups of $A$.

\bigskip

 We can now explain the construction of 
  Connes-Karoubi character \cite{MR972606}, see also \cite[Sec. 7.3]{MR2409415}.
The algebra $\cM^{d}$ has a natural Fr\'echet structure so that the notions of topological and relative $K$-theory used in   \cite{MR972606}  or \cite{MR2409415} coincide with our versions, see Remark \ref{kdjqkldjqlwkdjwqdwqdqwdwqdwd}.
  We start with the diagram derived from the right part of \eqref{jdklqwjdlkjqwdwd} and  the upper sequence in \eqref{slkkjle} (see \cite[4.10]{MR972606})
\begin{equation}\label{ewfwfewfewfeeee}\xymatrix{K^{top}_{d+2}(\cM^{d})\ar@{..>}@/_1.5cm/[dd]\ar[d]^{\ch_{top}^{cont}}\ar[r]&K_{d+1}^{rel}(\cM^{d})\ar[d]^{\ch^{cont}_{rel}}\ar[r]^{\alpha}&K_{d+1}(\cM^{d})\ar[r]^{0}\ar@{..>}[dd]^{\delta}&K_{d+1}^{top}(\cM^{d})\\HC^{cont}_{d+2}(\cM^{d})\ar[r]^{S}&HC^{cont}_{d}(\cM^{d})\ar[d]^{\phi_{d}}& & \\ \Z\ar[r]&\C\ar[r]&\C/ \Z&}\ .\end{equation}
It is a theorem of Karoubi \cite{MR816236} that the right upper map (marked  by $0$) vanishes. 
The map $\phi_{d}$ is given by the pairing with an explicit continuous cocycle
$\phi_{d}\in HC_{cont}^{d}(\cM^{d})$ which we will describe in  \eqref{dwedwedssssssss} below.
It has been verified in   \cite{MR972606} that elements coming from $\Z\cong K^{top}_{d+2}(\cM^{d})$
are mapped to integers under the obvious composition indicated by the left dotted arrow.
The right dotted arrow  is the multiplicative character. It is defined by the obvious diagram chase.

\begin{rem}\label{ijwdiwjqdlqwdwqdwqdwq}{\rm
In this remark we describe the cocycle $\phi_{d}$ explicitly.
The formula will be used in  the Remarks \ref{kjclkcjdscdscsdcdscdc} and \ref{dhedwedwedwedwedewdewd} below. Our description of  $\phi_{d}$
 employs the chain complex
$C^{\lambda,cont}_{*}(A)$  given in \cite[2.1.4]{MR1600246}
in order to  calculate $HC^{cont}_{*}(A)$ for a unital    locally convex algebra over $\C$. In particular
$C^{\lambda,cont}_{n}(A)$ is the space of coinvariants for the action of the cyclic permutation group on $A^{\otimes_{\pi} n+1}$. We use the notation $[a^{0}\otimes\dots\otimes a^{n}]$ in order to denote elements in $C^{\lambda,cont}_{n}(A)$. 

Furthermore, for a locally convex algebra $A$ we calculate the cyclic cohomology
$HC_{cont}^{d}(A)$ using the complex $C_{\lambda,cont}^{*}(A)$, where
$C_{\lambda,cont}^{n}(A)$ is the $\C$-vector space of continuous
 multilinear and cyclically invariant maps $A^{\times n+1}\to \C$.
We have a natural pairing
$$C_{\lambda,cont}^{n}(A)\times C^{\lambda,cont}_{n}(A)\to \C$$
given by
$$(\phi,[a^{0}\otimes\dots\otimes a^{n}])\to \phi(a^{0},\dots,a^{n})\ .$$
Using these conventions the map $\phi_{d}:HC^{cont}_{d}(\cM^{d})\to \C$ in \eqref{ewfwfewfewfeeee} is given  for odd $d$ by the pairing with the cocycle
 (using the notation introduced in Remark \ref{kljkljddlqwdqwdqwd})
\begin{equation}\label{dwedwedssssssss}\phi_{d}(a^{0},\dots,a^{d}):=(-1)^{\frac{d-1}{2}} \frac{d!}{(2\pi i)^{\frac{d-1}{2}}(\frac{d-1}{2})!} \Tr\left[ z \left(\begin{array}{cc}0&a^{0}_{12}\\a_{21}^{0}&0\end{array}\right) \dots\left(\begin{array}{cc}0&a^{d}_{12}\\ a^{d}_{21}&0\end{array}\right) \right]\ ,\end{equation}
where
$$z:= \left(\begin{array}{cc}1&0\\0&-1\end{array}\right)\ .$$
 }\hB
\end{rem}

 \begin{rem}\label{kjclkcjdscdscsdcdscdc}
{\rm   In this remark we approach an explicit formula for the composition
$$\delta\circ K(b_{\Dirac})\circ \partial:K^{rel}_{d+1}(C^{\infty}(X))\to  \C/\Z\ .$$  
Here for a homomorphism $b$ between algebras  we denote by $K(b)$ or $HC(b)$  the induced maps in $K$-theory or cyclic homology.
In view of  \eqref{ewfwfewfewfeeee} and the naturality of $\ch^{cont}_{rel}$ we have the equality
$$\delta\circ K(b_{\Dirac})\circ  \partial=[-]_{\C/\Z}\circ \phi_{d}\circ \ch^{cont}_{rel}\circ K(b_{\Dirac})=[-]_{\C/\Z}\circ \phi_{d}\circ HC(b_{\Dirac})\circ \ch^{cont}_{rel}\ .$$

It is clearly complicated  to write down an explicit formula for the relative Chern character $\ch^{cont}_{rel}$.
But we can
  give an explicit formula for the composition  $$\phi_{d}\circ HC(b_{\Dirac}):
HC_{d}(C^{\infty}(X))\to HC_{d}(\cM^{d})\to \C\ .$$
  We continue with the notation introduced in Remark \ref{kljkljddlqwdqwdqwd}.
For $f\in C^{\infty}(X)$ we have
$$[(P^{+}-P^{-}),f]=2 (P^{+}fP^{-}-P^{-}fP^{+})=2\left(\begin{array}{cc}0&P^{+}fP^{-}\\-P^{-}fP^{+}&0\end{array}\right)\ .$$
Combining \eqref{dwedwedssssssss} with \eqref{dqwdwqdwqdqwdqwdw324324234} we see that  $\phi_{d}\circ HC(b_{\Dirac})$ is represented by the cochain 
\begin{equation}\label{smwnqsmwqnsqws} (f_{0},\dots,f_{d})\mapsto  
   \frac{-2^{d+1}  d!}{(2\pi i)^{\frac{d-1}{2}}(\frac{d-1}{2})!}
\Tr\left((P^{+}-P^{-})[(P^{+}-P^{-}),f_{0}]\dots [(P^{+}-P^{-}),f_{d}]\right)\end{equation}
This formula is a first step in the direction of the main result of \cite{MR2817643}.
On the other hand it is still a complicated non-local formula. By  standard methods of local index theory using 
 e.g. the heat kernel and Getzler rescaling one can produce local  cocyles  representing the same cohomology class, see e.g. \cite{MR1047274}, \cite{MR1334867}.
In Lemma \ref{djqwjnkddqwdqwd} we avoid complicated analysis and the 
struggle with normalizations by using the Atiyah-Singer index theorem.
The explicit local formula will be stated in Remark \ref{dhedwedwedwedwedewdewd}.

%

}\hB
\end{rem}

\subsection{The conjecture}

The connective covering morphism of spectra \eqref{egdjhdgjhd} induces a morphism of spectra
\begin{equation}\label{hdghasgdjasdsdhasd87987}\ku\C/\Z\to \bku\C/\Z\ .\end{equation}

%
%
%

 Let $X$ be a   closed Riemannian manifold  of odd dimension $d$ equipped with a generalized Dirac
operator $\Dirac$. Then we define the map $r_{\Dirac}$ as the following composition:
 \begin{equation}\label{wdqwdwqdwqdwqdq32432432434}r_{\Dirac}:\ku\C/\Z^{-d-2}(X)\stackrel{\eqref{hdghasgdjasdsdhasd87987}}{\to} \bku\C/\Z^{-d-2}(X)\stackrel{\iota_{d+1}}{\cong} \bku\C/\Z^{-1}(X) \cong 
  \widehat{KU}^{0}_{flat}(X)\stackrel{ \rho_{\Dirac}}{\to} \C/\Z\ ,\end{equation}
  where $\iota_{2k}:\bku\C/\Z^{p}(X)\stackrel{\cong}{\to} \bku\C/\Z^{p+2k}(X)$ is again a periodicity operator.
We  now consider the   diagram:
 \begin{equation}\label{hdjkwhkdqwdqwdqwdwq62876348}\xymatrix{& \ku\C/\Z^{-d-2}(X)\ar[dr]^{r_{\Dirac}}&\\K_{d+1}(C^{\infty}(X))\ar[ur]^{\sigma_{d+1}}\ar[dr]^{b_{\Dirac}}&&\C/\Z\\& 
 K_{d+1}(\cM^{d})\ar[ur]^{\delta}&}\ ,\end{equation}
 where   $b_{\Dirac}$ is defined in \eqref{fevgerhjvhjbjhervvev} and classifies the Fredholm module of $\Dirac$, $\delta$ is the multiplicative character of Connes-Karoubi \eqref{edwededwedwe32423423554345}, and $\sigma_{d+1}$ is defined in Definition \ref{kdmlkqmdlqwdqwdwqd}.


 \begin{con}\label{lkwdjwqlkdqwldwqdwqww}   Let $X$ be a   closed Riemannian manifold  of odd dimension $d$ equipped with a generalized Dirac
operator $\Dirac$. Then the  diagram \eqref{hdjkwhkdqwdqwdqwdwq62876348} commutes.
 \end{con}
    
 In the present paper we show this conjecture for topologically trivial classes in $K_{d+1}(C^{\infty}(X))$. The precise formulation of this result is Theorem \ref{dlkdjqwkdjwqlkdqwdqwd}.   

\subsection{Comparison of certain cocycles}\label{dwlkdwqldqwdqwdqwdwdwd}

In this subsection we prepare the proof of Theorem \ref{dlkdjqwkdjwqlkdqwdqwd} by providing a differential geometric formula for the composition $\delta\circ HC(b_{\Dirac})$. The main result of this subsection is Lemma \ref{djqwjnkddqwdqwd}.

 In the following we define the cyclic homology  $HC_{*}(A)$  of an  associative unital  algebra $A$ over $\C$      as the homology of the standard cyclic complex $CC_{*}(A)$. For details we refer to \cite[2.1.9]{MR1600246}  where this complex is denoted by $\Tot \:\cB(A)$.
 Explicitly, we define
$$CC_{n}(A):= \left\{ \begin{array}{cc}\bigoplus_{k=0}^{\frac{n}{2}} A^{\otimes 2k+1}& \mbox{$n$ even}\\[0.3cm]
\bigoplus_{k=0}^{\frac{n-1}{2}} A^{\otimes 2k}&\mbox{$n$ odd}\end{array}\right. \ . $$
As in Subsection \ref{lkewdjlkwefwfwfewffefewfewf}, for a unital locally convex  algebra we define the continuous cyclic homology complex
$CC^{cont}(A)$ and its homology $HC^{cont}_{*}(A)$ similarly but using projective tensor products $\otimes_{\pi}$ instead of algebraic ones, see Remark \eqref{ldjlkqdjqwldqwdwqd7987}.

\begin{rem}{\em 
As shown in \cite[2.1.4]{MR1600246} there is    a natural quasi-isomorphism
\begin{equation}\label{chjkhdkhhkqjwhdwqd}CC(A)\to C^{\lambda}(A)\end{equation} which we use in order to compare the present definition of cyclic homology with the one used in  Remark \ref{ijwdiwjqdlqwdwqdwqdwq}. The  quasi-isomorphism \eqref{chjkhdkhhkqjwhdwqd} is induced by a chain-complex level  projection map, which in   degree $n$  is  the homomorphism
$CC_{n}(A)\to C^{\lambda}_{n}(A)$  given  by  (we write the formula for odd $n$)   
\begin{equation}\label{gdqjhgdjqwgdwqjdgjwqgdjqwdwqddw}\oplus_{k=0}^{\frac{n-1}{2}} \:\: a_{0}^{k}\otimes\dots\otimes a^{k}_{2k+1}\mapsto [a_{0}^{\frac{n-1}{2}}\otimes\dots\otimes a^{\frac{n-1}{2}}_{n}] \ .\end{equation}
There is a similar quasi-isomorphism in the continuous case.
 }\hB
\end{rem}

We define  a morphism of graded groups (see \eqref{mdmwwqdqwd87qw9} for the definition  of $DD(X)$)
\begin{equation}\label{hghdgj786786hjhd23d}\pi:CC^{cont}(C^{\infty}(X))\to  DD(X) \end{equation}
by the following prescription:
\begin{enumerate} \item If $n$ is odd, then we define
$CC^{cont}_{n}(C^{\infty}(X))\to \prod_{p\in \Z} (\sigma^{\le p}\Omega)[2p]^{-n}(X)$ by
\begin{equation}\label{gdqjhgdjqwgdwqjdgjwqgdjqwdwqddw1}\bigoplus_{k=0}^{\frac{n-1}{2}}f^{k}_{0}\otimes \dots \otimes f^{k}_{2k+1}\mapsto    
\sum_{k=0}^{\frac{n-1}{2}}%
 \frac{ b^{\frac{n+1}{2}+k}}{(2k+1)!}  f^{k}_{0}df^{k}_{1}\wedge \dots \wedge df^{k}_{2k+1}\ .\end{equation}
\item
If $n$ is even, then we define
$CC^{cont}_{n}(C^{\infty}(X))\to \prod_{p\in \Z} (\sigma^{\le p}\Omega)[2p]^{-n}(X)$ by
$$\bigoplus_{k=0}^{\frac{n}{2}}f^{k}_{0}\otimes \dots \otimes f^{k}_{2k}\mapsto   \sum_{k=0}^{\frac{n} {2}}  \frac{b^{\frac{n}{2}+k}}{(2k)!}  f^{k}_{0}df^{k}_{1}\wedge \dots \wedge df^{k}_{2k}\ .$$\\
\end{enumerate}
In these formulas we use the variable $b$ of degree $-2$ and the identification
$$\prod_{p\in \Z} \Omega[2p](X)\cong \Omega[b,b^{-1}](X)\ .$$
Under this identification the series
$\sum_{p\in \Z} b^{p}\omega(p)\in \Omega[b,b^{-1}](X)$ corresponds to the family
$(\omega(p))_{p\in \Z} \in \prod_{p\in \Z} \Omega[2p](X)$.
By  \cite[2.3.6]{MR1600246} the map $\pi$ is a morphism of chain complexes. By the calculation of the continuous  cyclic homology of $C^{\infty}(X)$ by Connes $\pi$  is actually a quasi-isomorphism.

\bigskip

We have a projection \begin{equation}\label{akljclkasjcasc}\psi:DD^{per}(X)\to DD(X)\end{equation} induced by  the projections in the components
$$DD^{per}(p)\to DD(p)\ , \quad \Omega[2p]\to (\sigma^{\le p}\Omega)[2p]$$
for all $p\in \Z$.
\begin{rem}{\rm
This projection \eqref{akljclkasjcasc} must not be confused with the projection \eqref{jkldjwqdlkwqjdqwdwdqwd}.
The latter is given by $$DD^{per}(X)\stackrel{\eqref{akljclkasjcasc}}{\to} DD(X)\stackrel{S}{\to} DD(X)[2]\ ,$$ where
$S((\omega(p))_{p\in \Z})=(\omega(p+1))_{p\in \Z}$.
}\hB\end{rem}

We now use that  
 $d\in \nat$ is odd and that  $\dim(X)=d$. Under these assumptions the map $\psi$ in 
  \eqref{akljclkasjcasc} induces   an isomorphism \begin{equation}\label{r4r34r34r43r4343343} \psi_{d}:HP^{-d}(X)\stackrel{\mbox{\tiny def}}{=}  H^{-d}(DD^{per}(X))\stackrel{\psi}{\to} H^{-d}(DD(X)) \ .\end{equation} We consider the isomorphism $\pi_{d}$ defined as the following composition of isomorphisms
  \begin{equation}\label{pid}\pi_{d}: HC^{cont}_{d}(C^{\infty}(X))\stackrel{\pi }{\longrightarrow} H^{-d}(DD(X)) \stackrel{\psi_{d}^{-1}}{\to}    HP^{-d}(X) \stackrel{\iota_{d-1}  }{\longrightarrow} 
HP^{-1}(X)\ .\end{equation}

\bigskip 


%

 Let $\Dirac$ be a generalized Dirac operator on $X$. Using the local index density $\hA(\Dirac)$ (see Remark \ref{ewjflewjflewfewfewf})
we define the map
\begin{equation}\label{kklwdjlkqwdqwdwdwqdwqd13}\tilde \rho_{\Dirac}:HP^{-1}(X)\to \C\ , \quad \tilde \rho_{\Dirac}([\gamma]):=\int_{X} \hA(\Dirac)\wedge \iota_{d+1}\gamma \ .\end{equation}

\begin{lem}\label{djqwjnkddqwdqwd}
Let $X$ be a closed manifold of odd dimension $d$ and $\Dirac$ be a generalized Dirac operator on $X$. Then the square
$$\xymatrix{HC^{cont}_{d}(C^{\infty}(X))\ar[r]^{  HC(b_{\Dirac})}\ar[d]^{\pi_{d}}&HC_{d}^{cont}(\cM^{d})\ar[d]^{\phi_{d}}  \\
HP^{-1}(X)\ar[r]^{\tilde \rho_{\Dirac} }&\C }$$
commutes.
\end{lem}
\proof

 Our task is  to compare the composition of the quasi-isomorphism \eqref{chjkhdkhhkqjwhdwqd} with the map \eqref{smwnqsmwqnsqws} on the one hand, and  the map $\tilde \rho_{\Dirac}$ defined in \eqref{kklwdjlkqwdqwdwdwqdwqd13} on the other.
It seems to be difficult to do this by an explicit calculation. Therefore we give an indirect argument based on the Atiyah-Singer index theorem. Our argument will not use explicit formulas. The convention for fixing normalizations described in Remark \ref{ewjflewjflewfewfewf} automatically takes care of the correct normalizations of $\hA(\Dirac)$ and $\phi_{d}$ in Remark \ref{ijwdiwjqdlqwdwqdwqdwq}.

\bigskip

We consider the composition of the map marked by $!!!$ in \eqref{jkwhdkqwdqwdwqd} with the Chern character   given by the third column in \eqref{slkkjle} in the case $A=\C$:
$$\bK^{top}_{C^{\infty}(X)}\stackrel{!!!}{\to} i_{X}^{*}\bK^{top}_{\C}\stackrel{\ch^{cont}_{per}}{\to} i_{X}^{*}\CC^{cont,per}_{\C}\ .$$
By evaluation at $*$ and taking the $(d+2)$'th homotopy group
we obtain the left triangle in the following diagram:
\begin{equation}\label{kjdklqwjdlkqwdqwd}\xymatrix{K_{d+2}^{top}(C^{\infty}(X)) \ar[rdd]\ar[r]&\pi_{d+2}(\bK^{top}_{\C}(X))\ar@{..>}[ddr]^(0.3){\iota_{d+1}\circ \ch^{cw}_{per}}\ar[r]^{!}\ar[dd]&HC^{cont}_{d}(C^{\infty}(X))\ar[r]^{  HC(b_{\Dirac})  }\ar[dd]^{\pi_{d}}& HC^{cont}_{d}(\cM^{d})\ar[dd]^{\phi_{d}}\\&&& \\&HC^{cont,per}_{d+2}(C^{\infty}(X))\ar[uur]^(0.3){q}\ar[r]&
HP^{-1}(X)\ar[r]^{\tilde\rho_{\Dirac} }&\C }\ .\end{equation}
The map marked by $q$ is induced by the map marked by this symbol in \eqref{llwdlqwddqwdkllkklklklkll}.
By construction the three solid  triangles on the left of \eqref{kjdklqwjdlkqwdqwd} commute.

Let  $K^{C^{*}}(-)$ denote the  usual $K$-theory for  $C^{*}$-algebras \cite{MR1656031}.
Since $X$ is $d$-dimensional, the connective covering \eqref{egdjhdgjhd} induces an isomorphism marked by $*$ in the following chain of isomorphisms:
\begin{equation}\label{dedwdwedwed8789733333}\pi_{d+2}(\bK_{\C}^{top}(X))\cong \ku^{-d-2}(X)\stackrel{*}{\cong} \bku^{-d-2}(X)
\cong K^{C^{*}}_{d+2}(C(X))
\ .\end{equation}Under this identification the  map
$$  \bku^{-d-2}(X)\to HP^{-1}(X)$$  induced by the dotted arrow in \eqref{kjdklqwjdlkqwdqwd} is the composition   $\iota_{d+1}\circ \ch^{cw}_{per}$  of the  usual Chern character and the shift, use Lemma \ref{jjljldkqjwdqwdqwd}.
In particular, its image is a lattice of full rank in $HP^{-1}(X)$.
Since $\pi_{d}$ is an isomorphism, in order to show that the right square in \eqref{kjdklqwjdlkqwdqwd} commutes is suffices to verify that the right hexagon (omit the left upper corner) commutes.

We consider the map
\begin{equation}\label{wjdqjlkwjlqkjdqwd87897}  K_{d+2}^{C^{*}}(C(X))\stackrel{\eqref{dedwdwedwed8789733333}}{\cong} \pi_{d+2}(\bK^{top}_{\C}(X)) \stackrel{!}{\to} HC^{cont}_{d}(C^{\infty}(X)) \end{equation} in the upper line of  \eqref{kjdklqwjdlkqwdqwd}.
The cocycle $\phi_{d}$ is normalized exactly such that the composition
of \eqref{wjdqjlkwjlqkjdqwd87897} with $\phi_{d}\circ HC(b_{\Dirac})$ is  the integer-valued function obtained by the pairing of $K_{d+2}^{C^{*}}(C(X))$
with the Fredholm module of $\Dirac$.

%


%

The  down-right-composition in the right hexagon of \eqref{kjdklqwjdlkqwdqwd} maps the $K$-theory class $x\in K_{d+2}^{C^{*}}(C(X))$ to  \begin{equation}\label{hgwshjqgwsqws876qws876s8768s6wqs}\int_{X} [\hA(\Dirac)]\cup \iota_{d+1}\ch^{cw}_{per}(x)\end{equation}which is apriori a complex number.
 The Atiyah-Singer index theorem encoded in equation \eqref{kdklqjwdljwqdlwqdwqdqwdqwdwd} shows that \eqref{hgwshjqgwsqws876qws876s8768s6wqs} is an integer and equal to the index pairing.
 So the  down-right-composition coincides with  the  right-down-composition.
\hB

\begin{rem}\label{dhedwedwedwedwedewdewd}{\rm  This is a continuation of Remark \ref{kjclkcjdscdscsdcdscdc}.
The following two cocycles  (a) and (b) on $CC_{d}^{cont}(C^{\infty}(X))$ represent the same 
map $HC_{d}^{cont}(C^{\infty}(X))\to \C$. We describe result of the application of the two cocyles to the chain
$$\oplus_{k=0}^{\frac{d-1}{2}} \:\: f_{0}^{k}\otimes\dots\otimes f^{k}_{2k+1}\in CC_{d}^{cont}(C^{\infty}(X))\quad :$$
\begin{eqnarray*}
 (a)&& \frac{-2^{d+1}  d!}{(2\pi i)^{\frac{d-1}{2}}(\frac{d-1}{2})!}
\Tr\left((P^{+}-P^{-})[(P^{+}-P^{-}),f^\frac{d-1}{2}_{0}]\dots [(P^{+}-P^{-}),f^{\frac{d-1}{2}}_{d}]\right)\\(b)&&
\sum_{k=0}^{\frac{d-1}{2}}\frac{1}{(2k+1)!}\int_{X} [\hA(\Dirac)]_{d-2k-1} \wedge  f^{k}_{0}df_{1}^{k}\wedge \dots\wedge d f^k_{2k+1 }
\end{eqnarray*}
The formula (a)   is one for 
$\phi_{d}\circ CC(b_{\Dirac})$ obtained by combining \eqref{gdqjhgdjqwgdwqjdgjwqgdjqwdwqddw}
and \eqref{smwnqsmwqnsqws}. The  formula (b) gives   $\tilde \rho_{\Dirac}\circ \pi_{d}$ and is derived from a combination of \eqref{kklwdjlkqwdqwdwdwqdwqd13} and \eqref{gdqjhgdjqwgdwqjdgjwqgdjqwdwqddw1}. The equality of the cohomology classes of (a) and (b) is the assertion of Lemma \ref{djqwjnkddqwdqwd}.

}\hB
\end{rem}

 \subsection{Proof of the conjecture for topologically trivial classes}\label{djoiqwjdoqwdjwqdwqdd}

In this subsection we prove Theorem \ref{dlkdjqwkdjwqlkdqwdqwd}.
We must show the equality of homomorphisms
$$r_{\Dirac}\circ \sigma_{d+1}\circ \partial=\delta\circ K(b_{\Dirac})\circ \partial:K_{d+1}^{rel}(C^{\infty}(X))\to \C/\Z\ .$$
This goal is achieved by the following chain of equalities:
\begin{eqnarray}\lefteqn{
r_{\Dirac}\circ \sigma_{d+1}\circ \partial}&&\nonumber\\
&=&r_{\Dirac}\circ  a  \circ \reg^{rel}_{X}\label{hhwkjdhqwkhdqwdwqdwqd1}\\
&=&r_{\Dirac}\circ  a \circ \psi_{d}\circ \psi_{d}^{-1}  \circ \reg^{rel}_{X}\label{hhwkjdhqwkhdqwdwqdwqd2}\\
&=&r_{\Dirac}\circ  a \circ \psi_{d} \circ \iota_{-d+1}\circ \iota_{d-1}\circ \psi_{d}^{-1}  \circ \reg^{rel}_{X}\label{hhwkjdhqwkhdqwdwqdwqd3}\\
&=&\rho_{\Dirac} \circ  a_{\widehat{\bku}} \circ  \iota_{d-1}\circ \psi_{d}^{-1}  \circ \reg^{rel}_{X}\label{hhwkjdhqwkhdqwdwqdwqd4}\\
&=&[\cdots]_{\C/\Z}\circ \tilde\rho_{\Dirac} \circ \iota_{d-1}\circ \psi_{d}^{-1}  \circ \reg^{rel}_{X}\label{hhwkjdhqwkhdqwdwqdwqd5}\\
&=&[\cdots]_{\C/\Z}\circ \tilde\rho_{\Dirac} \circ  \iota_{d-1}\circ \psi^{-1} \circ \pi\circ \pi^{-1}  \circ \reg^{rel}_{X}\label{hhwkjdhqwkhdqwdwqdwqd6}\\&=&
[\cdots]_{\C/\Z}\circ \tilde\rho_{\Dirac} \circ  \pi_{d}\circ \pi^{-1}  \circ \reg^{rel}_{X}\label{hhwkjdhqwkhdqwdwqdwqd7}\\
&=&[\cdots]_{\C/\Z}\circ\phi_{d}\circ HC(b_{\Dirac})\circ \pi^{-1}  \circ \reg^{rel}_{X}\label{hhwkjdhqwkhdqwdwqdwqd8}\\
&=&[\cdots]_{\C/\Z}\circ\phi_{d}\circ HC(b_{\Dirac})\circ \pi^{-1}  \circ \pi\circ \ch^{cont}_{rel}\label{hhwkjdhqwkhdqwdwqdwqd9} \\
&=&[\cdots]_{\C/\Z}\circ\phi_{d}\circ HC(b_{\Dirac})\circ  \ch^{cont}_{rel} \label{hhwkjdhqwkhdqwdwqdwqd10}\\
&=&[\cdots]_{\C/\Z}\circ\phi_{d}\circ \ch^{cont}_{rel}\circ K^{rel}(b_{\Dirac})  \label{hhwkjdhqwkhdqwdwqdwqd11} \\
&=&\delta \circ K(b_{\Dirac})\circ \partial \label{hhwkjdhqwkhdqwdwqdwqd12}
 \end{eqnarray}
 In the following we explain the steps in detail:
 \begin{enumerate}
 \item For \eqref{hhwkjdhqwkhdqwdwqdwqd1} we use Corollary \ref{hqkdjhqwkjdhqwd} and the Definition \ref{kdmlkqmdlqwdqwdwqd} of $\sigma_{d+1}$ in terms of $\hat \reg_{X}$.
 \item At \eqref{hhwkjdhqwkhdqwdwqdwqd2} we insert $\id=\psi^{-1}_{d}\circ \psi_{d}$, where $\psi_{d}$ is defined in \eqref{r4r34r34r43r4343343}.
 \item At \eqref{hhwkjdhqwkhdqwdwqdwqd3} we insert $\iota_{-d+1}\circ \iota_{d-1}=\id$, where $\iota_{k}$ is the periodicity operator
 introduced in \eqref{shiftiso} for every $k\in 2\Z$. Note that $d-1$ is even.
 \item At \eqref{hhwkjdhqwkhdqwdwqdwqd4} we use the commutative diagram
 $$\xymatrix{HP^{-1}(X)\ar[d]^{\iota_{-d+1}}\ar[r]^{a_{\widehat{\bku}} }&\widehat{KU}^{0}_{flat}(X)\ar@/^0.8cm/@{..>}[rrr]^{\rho_{\Dirac}}\ar[r]^{\cong}&\bku\C/\Z^{-1}(X)\ar[r]_{\iota_{-d-1}}^{\cong} &\bku\C/\Z^{-d-2}(X)&\C/\Z\\HP^{-d}(X)\ar[r]^(0.5){\psi_{d}}& \pi_{d+1}(\Sigma \DD(X))  \ar[r]^{a}  &\widehat{ku}_{flat}^{ -d-1}(X) \ar[r]^{\cong }&\ku\C/\Z^{-d-2}(X)\ar[u]_{\eqref{hdghasgdjasdsdhasd87987}}\ar[ur]^{r_{\Dirac}}& }\ . $$ 
  \item At \eqref{hhwkjdhqwkhdqwdwqdwqd5} we use the definition \eqref{kklwdjlkqwdqwdwdwqdwqd13} of $\tilde \rho_{\Dirac}$ and the observation based on formula \eqref{edewdwede}  that  its composition with $[\cdots]_{\C/\Z}$ coincides with  the composition $\rho_{\Dirac}\circ a_{\widehat{\bku}}$.
  \item At \eqref{hhwkjdhqwkhdqwdwqdwqd6} we insert $\pi\circ \pi^{-1}=\id$, where the isomorphism $\pi$ is defined in \eqref{hghdgj786786hjhd23d}.
  \item At \eqref{hhwkjdhqwkhdqwdwqdwqd7} we insert the definition   \eqref{pid} of $\pi_{d}$.
  \item At \eqref{hhwkjdhqwkhdqwdwqdwqd8} we use Lemma \ref{djqwjnkddqwdqwd}.
  \item At \eqref{hhwkjdhqwkhdqwdwqdwqd9} we use the equality
$\reg_{X}^{rel}=\pi\circ \ch^{cont}_{rel}$.
\item At \eqref{hhwkjdhqwkhdqwdwqdwqd10} we delete $\pi^{-1}\circ \pi$.
\item At \eqref{hhwkjdhqwkhdqwdwqdwqd11} we use the naturality of $\ch^{cont}_{rel}$ expressed by the commutative diagram:
$$\xymatrix{K^{rel}_{d+1}(C^{\infty}(X))\ar[r]^{\ch_{rel}^{cont}}\ar[d]^{K^{rel}(b_{\Dirac})}&HC_{d}^{cont}(C^{\infty}(X))\ar[d]^{HC(b_{\Dirac})}\\ K^{rel}_{d+1}(\cM^{d})\ar[r]^{\ch^{cont}_{rel}}&HC^{cont}_{d}(\cM^{d})}\ .$$
\item Finally, at \eqref{hhwkjdhqwkhdqwdwqdwqd12} we use the definition of $\delta$ in terms of the commutative diagram \eqref{ewfwfewfewfeeee}
 \end{enumerate}
 \hB

\bibliographystyle{alpha}
\bibliography{smoothreg}
\end{document}